\documentclass[12pt]{amsart}
 \usepackage{amsmath}
 \usepackage{epsfig} 
 \usepackage{graphics} 

\def\a{\alpha}
\def\b{\beta}

\def\l{\langle}
\def\r{\rangle}
\def\Ga{\Gamma}
\def\la{\lambda}

\def\p{\partial}

\def\tr{\mbox{tr}}

\newtheorem{Theorem}{Theorem}[section]

\newtheorem{Lemma}{Lemma}[section]

\theoremstyle{definition}
\newtheorem{definition}{Definition}[section]
\newtheorem{Remark}{Remark}[section]

\catcode`@=11 \@addtoreset{equation}{chapter}

\catcode`@=12

\title[ ]
{\bf Estimates and Existence Results for a Fully Nonlinear Yamabe 
Problem on Manifolds with
Boundary}

\author[]{}

\thanks{}

\begin{document}
\maketitle

\centerline{\scshape Qinian Jin
\begin{footnote}
{Department of Mathematics, Rutgers University,
110 Frelinghuysen Rd.,
Piscataway, NJ 08854, qjin@math.rutgers.edu}
\end{footnote}, Aobing Li
\begin{footnote}
{ 
Department of Mathematics, University of Wisconsin,
480 Lincoln Drive,
Madison, WI 53706, aobingli@math.wisc.edu
}
\end{footnote} 
and YanYan Li
\begin{footnote}
{Department of Mathematics, Rutgers University,
110 Frelinghuysen Rd.,
Piscataway, NJ 08854, yyli@math.rutgers.edu.
Partially supported by NSF grant DMS-0401118}
\end{footnote}
}



\medskip

\section{\bf Introduction}
\setcounter{equation}{0}

Let $({\mathcal M}, g)$ be a smooth compact Riemannian manifold of
dimension $n\ge 3$. The Yamabe problem is to find
metrics conformal to $g$ with constant scalar curvature. The
problem has been solved through the work of Yamabe
\cite{Ya60}, Trudinger \cite{T68}, Aubin \cite{Au76} and Schoen
\cite{Sch84}. See \cite{LP87} for a survey.
The Yamabe and related problems have attracted much attention in
the last 30 years or so, see, e.g., \cite{SY94}, \cite{Au98} and
the references therein.  Analogues of the Yamabe problem on
compact manifolds $({\mathcal M}, g)$ with boundary have
been studied by Cherrier \cite{Ch84}, Escobar \cite{E92a, E, E96},
Han and Li \cite{HL, HL00}, Ambrosetti, Li and Malchiodi
\cite{ALM02} and others.

In recent years, fully nonlinear versions of the Yamabe problem have
received much attention since the work of Viaclovsky \cite{V00}. Let
$({\mathcal M}^n, g)$ be a smooth compact Riemannian manifold of
dimension $n\ge 3$. The Schouten tensor is defined as
$$
A_g:=\frac{1}{n-2}\left(Ric_g-\frac{1}{2(n-1)} R_g g\right),
$$
where $Ric_g$ and $R_g$ are the Ricci tensor and the scalar
curvature of $g$ respectively. Let $\la(A_g)=(\la_1(A_g), \cdots,
\la_n(A_g))$ denote the eigenvalues of $A_g$ with respect to $g$.
One interesting problem is to find conformal metrics on
$({\mathcal M}, g)$ with a prescribed symmetric function of the
eigenvalues of the Schouten tensors.

To be more precise, let
$$
\Ga_1:=\left\{\la=(\la_1, \cdots, \la_n)\in {\mathbb R}^n: \sum
\la_i>0\right\}
$$
and
$$
\Ga_n:=\left\{\la=(\la_1, \cdots \la_n)\in {\mathbb R}^n: \la_i>0
\mbox{ for } 1\le i\le n\right\}.
$$
We assume that
\begin{equation}\label{1.1}
\Ga\subset {\mathbb R}^n \mbox{ is an open convex symmetric cone with
vertex at the origin}
\end{equation}
satisfying
\begin{equation}\label{1.2}
\Ga_n\subset \Ga\subset \Ga_1
\end{equation}
and assume that
\begin{equation}\label{1.3}
f\in C^\infty(\Ga)\cap C^0({\overline{\Ga}}) \mbox{ is a symmetric
function, $f>0$ in $\Ga$}
\end{equation}
verifying some of the following properties which will be specified in
each situation:
\begin{equation}\label{1.4}
f \mbox{ is homogeneous of degree one on } \Ga,
\end{equation}
\begin{equation}\label{1.5}
f_i:=\frac{\p f}{\p \la_i}>0 \mbox{ on } \Ga,
\end{equation}
and
\begin{equation}\label{1.6}
f \mbox{ is concave on } \Ga.
\end{equation}

Conditions (\ref{1.4}), (\ref{1.5}) and (\ref{1.6}) imply 
that (see \cite[Lemma 3.2]{Urbas})
\begin{equation}\label{1.6.5}
{\mathcal T}:=\sum_i f_i(\la)\ge c_0>0, \quad \forall \la\in \Ga
\end{equation}
for some positive number $c_0>0$.

The fully nonlinear Yamabe problem on a closed manifold
$({\mathcal M}, g)$ is to find a metric $\tilde{g}$ conformal to
$g$ such that
\begin{equation}\label{1.7}
F(A_{\tilde{g}}):=f(\la(A_{\tilde{g}}))=1 \quad \mbox{and} \quad
\la(A_{\tilde{g}})\in \Ga \quad \mbox{on } {\mathcal M}.
\end{equation}
When $(f, \Ga)=(\sigma_k^{1/k}, \Ga_k)$, this problem is known as
the $\sigma_k$-Yamabe problem in the literature, where, for each
$1\le k\le n$, $\sigma_k$ is the $k$-th elementary symmetric
function defined by
$$
\sigma_k(\lambda):=\sum_{i_1<\cdots<i_k}
\lambda_{i_1}\cdots\lambda_{i_k} \quad \mbox{for all }
\lambda=(\lambda_1, \cdots, \lambda_n)\in {\mathbb R}^n
$$
and
$$
\Ga_k:=\{\lambda\in {\mathbb R}^n: \sigma_l(\lambda)>0, 1\le
l\le k\}
$$
which is an open convex symmetric cone with vertex at
the origin.

For problem (\ref{1.7}) a number of existence results have been available
in the literature. In \cite{Via02} Viaclovsky established the
existence result for (\ref{1.7}) with $(f, \Ga)=(\sigma_n^{1/n},
\Ga_n)$ for a class of manifolds. In
\cite{CGY02, CGY03}, Chang, Gursky and Yang obtained the existence
result on 4-manifolds for (\ref{1.7}) with $(f,
\Ga)=(\sigma_2^{1/2}, \Ga_2)$. For $(f, \Ga)=(\sigma_k^{1/k},
\Ga_k)$ with $k=3, 4$ on $4$-manifolds and with $k=2, 3$ on
$3$-manifolds that are not simply connected, the existence result was
established by Gursky and Viaclovsky in \cite{GV04}. 
When $({\mathcal M}, g)$ is locally conformally flat and $(f,
\Ga)=(\sigma_k^{1/k}, \Ga_k)$ for $1\le k\le n$,  Guan-Wang
\cite{GW03a} and Li-Li \cite{LL03} independently proved the
existence of solutions of (\ref{1.7}). Li-Li \cite{LL03,
LL05} also established a general result that (\ref{1.7}) is still solvable if
$({\mathcal M}, g)$ is locally conformally flat and if $(f, \Ga)$
satisfies (\ref{1.1})--(\ref{1.6}) with $f|_{\p \Ga}=0$. In
\cite{GW03} Guan and Wang proved local interior $C^1$ and $C^2$
estimates for solutions of (\ref{1.7}) with $(f,
\Ga)=(\sigma_k^{1/k}, \Ga_k)$, such estimates were also studied in
\cite{LL03, GW04, STW05} and was extended to a general class of
$(f, \Ga)$ in \cite{Ch05}. Using such local estimates and the
algebraic fact found in \cite{GVW03} that $\la(A_g)\in \Ga_k$ for
$k>\frac{n}{2}$ implies the positivity of the Ricci tensor, Gursky and
Viaclovsky \cite{GV05} solved (\ref{1.7}) on general
manifolds if $(f, \Ga)$ satisfies (\ref{1.1})--(\ref{1.6}) with
$f|_{\p \Ga}=0$ and if $\Ga\subset \Ga_k$ for some
$k>\frac{n}{2}$. In \cite{TW05} Trudinger and Wang proved a
Harnack inequality for the set of metrics $\tilde{g}$
conformal to $g$ with $\la(A_{\tilde{g}})\in \Ga_k$ for some
$k>\frac{n}{2}$ and gave a different proof of
the existence result in \cite{GV05}. For $(f, \Ga)=(\sigma_k^{1/k},
\Ga_k)$ with $k\le \frac{n}{2}$ on general manifolds, Sheng,
Trudinger and Wang \cite{STW05} established the existence result
under a variational structure condition which includes the
case $k=2$ and $n\ge 4$, while Ge and Wang \cite{GeW05} independently
obtained a proof for $k=2$ and $n>8$.

In the rest of this paper we will assume that $({\mathcal M}^n, g)$, $n\ge 3$,
is a smooth compact Riemannian manifold with nonempty smooth
boundary $\p {\mathcal M}$. For a given constant $c\in {\mathbb R}$,
we are interested in finding a metric $\tilde{g}$ conformal to $g$
such that
\begin{equation}\label{1.8}
\left\{\begin{array}{lll}
F(A_{\tilde{g}}):=f(\la(A_{\tilde{g}}))=1, \quad
\la(A_{\tilde{g}})\in \Ga \quad \mbox{on } {\mathcal M},\\
\\
h_{\tilde{g}}=c \quad \mbox{on } \p {\mathcal M},
\end{array}
\right.
\end{equation}
where $h_{\tilde{g}}$ denotes the mean curvature of $\p {\mathcal
M}$ with respect to the outer normal (A Euclidean ball has positive
boundary mean curvature). Note that when $(f,
\Ga)=(\sigma_1, \Ga_1)$, this is the Yamabe problem on a
compact manifold with boundary. Therefore (\ref{1.8}) is a
fully nonlinear version of the Yamabe problem with boundary.

This problem was proposed by Li and Li in \cite{LL04a, LL04b} in
which they considered the corresponding blow-up problem of
(\ref{1.8}) and obtained some Liouville type theorems and a
Harnack type inequality. These results indicate positively that it
should be possible to establish some existence results for
(\ref{1.8}) under suitable conditions.

Writing $\tilde{g}=e^{-2u} g$ for some smooth function $u$ on
$\mathcal{M}$. Using the transformation laws for the Schouten tensor
and mean curvature, (\ref{1.8}) is equivalent to the fully
nonlinear elliptic equations
\begin{equation}\label{1.9}
\left\{\begin{array}{lll} F(U):=f(\la_g(U))=e^{-2u}, \quad
\la_g(U)\in \Ga \quad \mbox{on } {\mathcal M},\\
\\
\frac{\p u}{\p \nu}=c e^{-u}-h_g \quad \mbox{on } \p {\mathcal M}
\end{array}
\right.
\end{equation}
with
\begin{equation}\label{1.10}
U=\nabla^2 u+ d u\otimes d u -\frac{1}{2}|\nabla u|_g^2 g+A_g,
\end{equation}
where $\la_g(U)$ denotes the eigenvalues of $U$ with respect to $g$,
$\nu$ is the unit inward normal vector field to $\p
\mathcal{M}$ in $({\mathcal M}, g)$ and $\nabla$ denotes the
Levi-Civita connection with respect to $g$.

Recall that the second fundamental form $\Pi$ of $\p {\mathcal M}$
with respect to $g$ is defined as
$$
\Pi(X, Y):=-g(\nabla_X \nu, Y), \quad \mbox{for any } X, Y\in T(\p
{\mathcal M}),
$$
where $T(\p {\mathcal M})$ denotes the tangent bundle over $\p
{\mathcal M}$. A point $x\in \p {\mathcal M}$ is called an umbilic
point in $({\mathcal M}, g)$ if
$$
\Pi(X, Y)=h_g(x) g(X, Y) \quad \mbox{for all } X, Y\in T_x(\p
{\mathcal M}).
$$
The boundary $\p {\mathcal M}$ is called umbilic if every point of
$\p {\mathcal M}$ is an umbilic point. The notion of umbilic point
is conformally invariant, i.e. a point is umbilic with respect to
$g$ is still umbilic with respect to the metric
$\tilde{g}:=e^{-2u}g$ for any function $u\in C^2(\mathcal{M})$
(see \cite{E}).

Our first existence result is for locally conformally flat
manifolds with umbilic boundary.

\begin{Theorem}\label{T1.1}
Assume that $(f, \Ga)$ satisfies (\ref{1.1})--(\ref{1.6}) with
$f|_{\p \Ga}=0$ and that $({\mathcal M}, g)$ is a smooth compact locally
conformally flat Riemannian manifold with smooth umbilic boundary 
$\p {\mathcal M}$. Suppose that $\la(A_g)\in \Ga$ on ${\mathcal M}$
and $h_g\ge 0$ on $\p {\mathcal M}$. Then problem (\ref{1.9}) with $c=0$ has
a solution $u\in C^\infty({\mathcal M})$.
\end{Theorem}

The existence of solutions of (\ref{1.7}) with $(f, \Ga)$
satisfying (\ref{1.1})--(\ref{1.6}) and $f|_{\p \Ga}=0$ has been
proved in \cite{LL03b, LL05} on compact locally conformally flat
manifolds without boundary. The proof of Theorem \ref{T1.1} is
based on \cite{LL03b, LL05}. By making use of the double of a
compact manifold, the problem in Theorem \ref{T1.1} reduces to a
corresponding problem on compact locally conformally flat manifold
without boundary. In order to establish $C^0$ estimates, via the
Harnack inequality obtained in \cite{LL03b}, we need to assume 
$h_g\ge 0$ on $\p {\mathcal M}$.

Recall that the Yamabe problem on compact manifolds $({\mathcal
M}, g)$ with boundary $\partial {\mathcal M}$ is to find a
conformally related metric $\tilde{g}$ of constant scalar
curvature on ${\mathcal M}$ and constant mean curvature on $\p
{\mathcal M}$. By writing $\tilde{g}=u^{\frac{4}{n-2}}g$ for some
positive smooth function on ${\mathcal M}$, this problem is
equivalent to finding a smooth positive solution $u$ to the
boundary value problem
\begin{equation}\label{1.9.4}
\left\{\begin{array}{lll} -\frac{4(n-1)}{n-2} \Delta_g u +R_g u=
a u^{\frac{n+2}{n-2}} \quad \mbox{on } {\mathcal M},\\
\\
-\frac{2}{n-2}\frac{\p u }{\p \nu} + h_g u= c u^{\frac{n}{n-2}}
\quad \mbox{on } \p {\mathcal M},
\end{array}\right.
\end{equation}
where $a$ and $c$ are constants. For any $a>0$ and any $c$, the existence of a
solution of (\ref{1.9.4}) has been proved in \cite{HL, HL00} under
the assumption that $({\mathcal M}, g)$ is of positive type and
satisfies one of the following assumptions:

\begin{enumerate}
\item[(i)] $({\mathcal M}^n, g)$, $n\ge 3$, is locally conformally
flat with umbilic boundary;

\item[(ii)] $n\ge 5$ and ${\mathcal M}$ is not umbilical.
\end{enumerate}
For any $a>0$, it was proved earlier in \cite{E, E96} that (\ref{1.9.4})
is solvable for $c=0$ and for at least one $c_+(a)>0$
and one $c_{-}(a)<0$ under the same assumptions. Here we call a 
manifold $({\mathcal M}, g)$ of positive type if
$$
\la_1({\mathcal M}):=\min_{\varphi\in H^1({\mathcal M})\backslash
\{0\}} \frac{\int_{\mathcal M}\left(|\nabla\varphi|_g^2+c(n)R_g
\varphi^2\right)+\frac{n-2}{2}\int_{\p {\mathcal M}} h_g \varphi^2
}{\int_{\mathcal M} \varphi^2}>0,
$$
where $c(n)=\frac{n-2}{4(n-1)}$.

An interesting question for (\ref{1.8}) is to identify good
conditions which guarantee the existence of a solution. 
Theorem \ref{T1.1} is  such an attempt, and it
shows that $h_g\ge 0$ on $\p {\mathcal M}$ is a sufficient
condition. Unlike the Yamabe problem with boundary, we tend to
believe that the hypothesis ``$h_g\ge 0$ on $\p {\mathcal M}$'' in 
Theorem \ref{T1.1} can not be replaced by ``$\la_1({\mathcal M})>0$''.

Our next result concerns (\ref{1.9}) with $c> 0$. We consider more 
general equation
\begin{equation}\label{5.2}
\left\{\begin{array}{lll} F(U):=f(\la_g(U))=\varphi_0 e^{-2u},
\quad
\la_g(U)\in \Ga \quad \mbox{on } {\mathcal M},\\
\\
\frac{\p u}{\p \nu}=h_0 e^{-u}-h_g \quad \mbox{on } \p {\mathcal
M},
\end{array}\right.
\end{equation}
where $\varphi_0\in C^\infty({\mathcal M})$
and $h_0\in C^\infty(\p{\mathcal M})$ are 
positive functions. This problem is equivalent
to finding a metric $\tilde{g}$ conformal to $g$ such that
$f(\la(A_{\tilde{g}})=\varphi_0$ on ${\mathcal M}$ and
$h_{\tilde{g}}=h_0$ on $\p {\mathcal M}$.

\begin{Theorem}\label{T1.2}
Assume that $(f, \Ga)$ satisfies (\ref{1.1})--(\ref{1.6}) with
$f|_{\p \Ga}=0$ and $\Ga\subset \Ga_k$ for some $k>\frac{n}{2}$.
Let $({\mathcal M}, g)$ be a smooth compact Riemannian manifold
with smooth boundary $\p {\mathcal M}$. Suppose that $\la(A_g)\in
\Ga$ on ${\mathcal M}$, $h_g\ge 0$ on $\p {\mathcal M}$, 
$\p {\mathcal M}$ is umbilic, and $({\mathcal M}, g)$ is locally conformally flat
near $\p {\mathcal M}$. Then for any
positive functions $\varphi_0\in C^\infty({\mathcal M})$ and
$h_0\in C^\infty(\p{\mathcal M})$ problem
(\ref{5.2}) has a solution $u\in C^\infty({\mathcal M})$.
\end{Theorem}

In \cite{GV05, TW05} more general equations than (\ref{1.7}) on
general closed manifolds have been solved when $\Ga\subset \Ga_k$
for some $k<\frac{n}{2}$. By using the double of a manifold we
adapt the Harnack inequality of Trudinger-Wang \cite{TW05} to
our situation. This, together with the $C^1$ and $C^2$ estimates
in sections 2-4, allows us to obtain the existence result by modifying
the degree argument in \cite{TW05}.

In sections 2-3 we establish under suitable conditions
on $(f, \Ga)$ some local $C^1$ and $C^2$ estimates for solutions
of the following more general equation
\begin{equation}\label{1.11}
\left\{\begin{array}{lll} F(U):=f(\la_g(U))=\psi(x,u), \quad
\la_g(U)\in \Ga \quad \mbox{on } {\mathcal O}_1,\\
\\
\frac{\p u}{\p \nu}=\eta(x, u)-h_g \quad \mbox{on } {\mathcal
O}_1\cap \p {\mathcal M},
\end{array}\right.
\end{equation}
where ${\mathcal O}_1$ is an open set of ${\mathcal M}$, $U$ is
defined by (\ref{1.10}), $\psi\in C^2({\mathcal O}_1\times
{\mathbb R})$ and $\eta\in C^2(({\mathcal O}_1\cap \p {\mathcal
M}) \times {\mathbb R})$. By extension we can
always assume that $\eta\in C^2({\mathcal O}_1\times {\mathbb
R})$.

$C^1$ and $C^2$ estimates have been studied extensively on
closed manifolds, see \cite{GW03, LL03, GW04, STW05, Ch05} for local
interior estimates and \cite{Via02} for global estimates. Global 
estimates have also been studied in \cite{G05} on
compact manifolds under Dirichlet boundary condition.

The first result on gradient estimates is the following.

\begin{Theorem}\label{T1.3}
Assume that $(f, \Ga)$ satisfies (\ref{1.1})--(\ref{1.5}),
(\ref{1.6.5}) and that $(\mathcal{M}, g)$ is a smooth compact Riemannian
manifold with smooth boundary $\p {\mathcal M}$. Let ${\mathcal O}_1$ be
an open set of ${\mathcal M}$ and let $u\in C^3(\mathcal{O}_1)$ be
a solution of (\ref{1.11}). If
\begin{equation}\label{1.12}
a\le u\le b \quad \mbox{on } {\mathcal O}_1
\end{equation}
for some constants $a$ and $b$, then, for any open set ${\mathcal
O}_2$ of ${\mathcal M}$ satisfying $\overline{\mathcal O}_2\subset
{\mathcal O}_1$,
$$
|\nabla u|_g\le C \quad \mbox{on } {\mathcal O}_2
$$
for some positive constant $C$ depending only on $n$ $(f, \Ga)$,
$g$, $\psi$, $\eta$, $a$, $b$, ${\mathcal O}_1$ and ${\mathcal
O}_2$.
\end{Theorem}

Note that the gradient estimate given in Theorem \ref{T1.3}
depends on the bound of $|u|$ on ${\mathcal O}_1$. If we only know
$u\ge -C_0$ on ${\mathcal O}_1$, the next result gives the gradient 
estimates for solutions of the equation
\begin{equation}\label{1.9.5}
\left\{\begin{array}{lll} F(U):=f(\la_g(U))=e^{-2u}, \quad
\la_g(U)\in \Ga \quad \mbox{on } {\mathcal O}_1,\\
\\
\frac{\p u}{\p \nu}=c e^{-u}-h_g \quad \mbox{on } {\mathcal
O}_1\cap \p {\mathcal M},
\end{array}
\right.
\end{equation}
where ${\mathcal O}_1$ is an open set of ${\mathcal M}$,
if $(f, \Ga)$ further satisfies the condition $(H_\a)$ introduced in
\cite{LL03}, see Definition \ref{D2.1} in section 2,

\begin{Theorem}\label{T2.1}
Assume that $(f, \Ga)$ satisfies (\ref{1.1})--(\ref{1.5}),
(\ref{1.6.5}) and the condition $(H_1)$, and that
$(\mathcal{M}, g)$ is a smooth compact Riemannian manifold with smooth
boundary $\p {\mathcal M}$. Let ${\mathcal O}_1$ be an open set of ${\mathcal
M}$ and let $u\in C^3({\mathcal O}_1)$ be a solution of
(\ref{1.9.5}). If
$$
u\ge -C_0 \quad \mbox{on } {\mathcal O}_1
$$
for some constant $C_0$, then, for any open set ${\mathcal O}_2$
of ${\mathcal M}$ satisfying $\overline{{\mathcal O}}_2\subset
{\mathcal O}_1$,
$$
|\nabla u|_g\le C \quad \mbox{on } {\mathcal O}_2
$$
for some positive constant $C$ depending only on $n$, $c$, $(f,
\Ga)$, $g$, $C_0$, ${\mathcal O}_1$ and ${\mathcal O}_2$.
\end{Theorem}

By assuming that $\p {\mathcal M}$ is umbilic and that $({\mathcal
M}, g)$ is locally conformally flat near $\p {\mathcal M}$, we derive
in section 3 we derive the following $C^2$ estimates for 
solutions of (\ref{1.11}) under suitable conditions on $\psi$ and
$\eta$.

\begin{Theorem}\label{T1.5}
Assume that $(f, \Ga)$ satisfies (\ref{1.1})--(\ref{1.6}) and that
$(\mathcal{M}, g)$ is a smooth compact Riemannian manifold with smooth
boundary $\p \mathcal{M}$. Suppose that $\p {\mathcal M}$ is umbilic and
$({\mathcal M}, g)$ is locally conformally flat near $\p {\mathcal M}$. 
Let ${\mathcal O}_1$ be an open set of
${\mathcal M}$ and let $u\in C^4({\mathcal O}_1)$ be a solution of
(\ref{1.11}). Assume that $\psi$ and $\eta$ satisfy one of the
following conditions:
\begin{enumerate}
\item[(i)] $\eta\equiv 0$, $\frac{\p \psi}{\p \nu}\equiv 0$ on
$({\mathcal O}_1\cap \p {\mathcal M})\times {\Bbb R}$ and $\psi\in
C^3({\mathcal O}_1\times {\Bbb R})$;

\item[(ii)] $\eta$ is positive on $({\mathcal O}_1\times \p
{\mathcal M})\times {\Bbb R}$ and $\psi$ is any function on
${\mathcal O}_1\times {\Bbb R}$.
\end{enumerate}
If
\begin{equation}\label{3.1}
|u|\le C_0 \quad \mbox{and}\quad |\nabla u|_g\le C_0 \quad
\mbox{on } {\mathcal O}_1
\end{equation}
for some constant $C_0$, then, for any open set ${\mathcal O}_2$
of ${\mathcal M}$ satisfying $\overline{\mathcal O}_2\subset
{\mathcal O}_1$,
$$
|\nabla^2 u|_g \le C \quad \mbox{on } {\mathcal O}_2
$$
for some constant $C$ depending only on $n$, $C_0$, $g$, $(f,
\Ga)$, $\psi$,  $\eta$, ${\mathcal O}_1$ and ${\mathcal O}_2$.
\end{Theorem}

Both Theorem \ref{T1.3} and Theorem \ref{T1.5} are used in the
proof of Theorem \ref{T1.2}. In section 4 we establish $C^2$ estimates 
on general manifolds with umbilic boundary, without any locally 
conformally flat assumption, for the following 
Monge-Amp\'{e}re type problem
\begin{equation}\label{4.1}
\left\{\begin{array}{lll} \det(g^{-1}\cdot U)=e^{-2n u}, \quad
\la_g(U)\in \Ga_n \quad \mbox{on } {\mathcal O}_1,\\
\\
\frac{\p u}{\p \nu}=c e^{-u}-h_g \quad \mbox{on } {\mathcal
O}_1\cap \p {\mathcal M}.
\end{array}
\right.
\end{equation}
where $U$ is defined by (\ref{1.10}).

\begin{Theorem}\label{T4.1}
Assume that $(\mathcal{M}, g)$ is a smooth compact Riemannian manifold
with smooth umbilic boundary $\p {\mathcal M}$. Let ${\mathcal O}_1$ be
an open set of ${\mathcal M}$ and let $u\in C^4({\mathcal O}_1)$
be a solution of (\ref{4.1}) with $c>0$. If
\begin{equation}\label{4.2}
|u|\le C_0 \quad \mbox{and} \quad |\nabla u|_g\le C_0 \quad
\mbox{on } \mathcal{O}_1
\end{equation}
for some constant $C_0$, then, for any open set ${\mathcal O}_2$
of ${\mathcal M}$ satisfying $\overline{\mathcal O}_2\subset
{\mathcal O}_1$,
$$
|\nabla^2 u|_g\le C \quad \mbox{on } {\mathcal O}_2
$$
for some positive constant $C$ depending only on $n$, $c$, $g$,
$C_0$, ${\mathcal O}_1$ and ${\mathcal O}_2$.
\end{Theorem}

The Dirichlet problem for equation (\ref{1.7}) has been studied
in \cite{G05} for $(f, \Ga)=(\sigma_k^{1/k}, \Ga_k)$ and the existence
of solutions is established whenever there exists an admissible 
supersolution. A similar problem for
$(f, \Ga)=(\sigma_n^{1/n}, \Ga_n)$ was studied in \cite{S05}. 
The Neumann problem for Hessian equations
has been studied in \cite{LTU86, T87, U95, CN99, SS03}, 
most of the works are for Monge-Ampere equations. The results in
\cite{T87} and \cite{U95} concern, respectively, general Hessian
equations on Euclidean balls and on general domains in dimension two.

We draw readers' attention to some closely related  independent
work of Sophie Chen in \cite{Ch06}.

\section{\bf Gradient estimates}
\setcounter{equation}{0}

In this section we will prove Theorem \ref{T1.3} and Theorem
\ref{T2.1} concerning gradient estimates for solutions of
(\ref{1.11}) and (\ref{1.9.5}). We use the distance function 
$d_g(x, \p {\mathcal M})$ in $({\mathcal
M}, g)$ to the boundary $\p {\mathcal M}$. Clearly
there is a suitable small constant $\delta_0>0$ such that $d_g(x,
\p \mathcal{M})$ is smooth in $\{x\in \mathcal{M}: d_g(x, \p
\mathcal{M})\le 2\delta_0\}$. Moreover
$$
\frac{\p}{\p \nu}d_g(x, \p \mathcal{M})=1\quad \mbox{on  } \p
\mathcal{M}.
$$

We will fix a positive constant $C_1$ such that
\begin{equation}\label{2.1}
\Pi(X, X)\ge -C_1 g(X, X)\quad \mbox{for }  X\in T(\p
\mathcal{M}).
\end{equation}

It is well-known that we can always find a metric conformal to $g$
with vanishing mean curvature on $\p {\mathcal M}$. Since a
conformal change of metrics does not affect our $C^1$ and $C^2$
estimates, without loss of generality, in sections 2--4 we 
always assume that $h_g=0$ on $\p {\mathcal M}$ in the arguments.

\begin{proof}[Proof of Theorem \ref{T1.3}]
By shrinking ${\mathcal O}_2$ if necessary, we can always choose a
cut-off function $\rho\in C_0^\infty({\mathcal O}_1)$ such that
\begin{equation}\label{2.1.1}
0\le \rho\le 1 \mbox{ in } {\mathcal O}_1, \quad |\nabla \rho|\le
C\sqrt{\rho} \mbox{ in } {\mathcal O}_1, \quad \rho=1 \mbox{ on }
{\mathcal O}_2
\end{equation}
and
\begin{equation}\label{2.1.2}
\frac{\p \rho}{\p \nu}=0 \quad \mbox{on } {\mathcal O}_1\cap \p
{\mathcal M} \mbox{ if } {\mathcal O}_1\cap \p {\mathcal M}\ne
\emptyset,
\end{equation}
where $C$ is a universal constant.

Let
$$
\gamma(t):=\frac{1}{\Lambda}(1+t-a)^{-\Lambda}, \quad t\in [a, b],
$$
where the number $\Lambda$ is large enough so that $\Lambda/(1+b-a)\ge 8$.
Then we choose a function $\varphi\in C^\infty({\mathcal M})$ such
that $\varphi(x)=d_g(x, \p {\mathcal M})$ when $d_g(x, \p
{\mathcal M})\le \varepsilon_0$, where $0<\varepsilon_0\le
\delta_0$ is sufficiently small. Since $u$ satisfies (\ref{1.12}),
we may further assume that $\varphi$ is chosen in a way so that
\begin{equation}\label{2.2}
|\eta_u \varphi|\le \frac{1}{2} \quad \mbox{and} \quad |\eta_{uu}
\varphi|\le \frac{1}{2} \quad \mbox{on } {\mathcal M}\times [a,
b].
\end{equation}

We now consider the function
$$
G:=\frac{1}{2}\rho \mbox{e}^{\a} |\nabla u-\b\nabla \varphi|_g^2
\quad \mbox{on } {\mathcal O}_1,
$$
where
$$
\a(x):=\la \varphi(x)+\gamma(u-\eta(x, u) \varphi), \quad
\b(x):=\eta(x, u)
$$
and $\la$ is a large number to be determined later. In order to
derive the desired bound on $|\nabla u|_g$ over ${\mathcal O}_2$,
it suffices to show that $G$ can be bounded in ${\mathcal O}_1$ by
some universal constant $C$ depending only on $n$, $(f, \Ga)$,
$g$, $\psi$, $\eta$, $a$, $b$, ${\mathcal O}_1$ and ${\mathcal
O}_2$. Suppose the maximum of $G$ over $\mathcal{O}_1$ is attained
at some point $x_0\in \mathcal{O}_1$. In the following we will
always assume that $G(x_0)\ge 1$; otherwise we are done.

We first claim that $x_0$ must be an interior point of
$\mathcal{O}_1$, i.e. $x_0\in \mathcal{O}_1\backslash\p
\mathcal{M}$. To this end, suppose $x_0\in \p \mathcal{M}$ and
choose an orthonormal frame field $\{e_1, \cdots, e_n\}$ around
$x_0$ such that $e_n=\nu$ on $\p \mathcal{M}$. In the following
for any smooth function $\phi$ we use $\phi_i, \phi_{ij}, \cdots$
to denote the covariant derivatives of $\phi$ of all orders. It is
easy to see that on $\p \mathcal{M}$ there hold
$$
\varphi_n=1, \quad \varphi_l=0 \quad \mbox{and} \quad
\varphi_{ln}=\varphi_{nl}=0 \quad  \mbox{for } 1\le l\le n-1.
$$
By using the boundary condition $u_n=\eta(x, u)$ we thus have on
$\p {\mathcal M}$ that
$$
u_n-\b \varphi_n=0 \quad \mbox{and} \quad
G=\frac{1}{2}\rho e^{\gamma(u)} |\nabla u-\b \nabla
\varphi|_g^2=\frac{1}{2}\rho e^{\gamma(u)} \sum_{l=1}^{n-1} u_l^2.
$$
Therefore, since $\rho_n=0$ on ${\mathcal O}_1\cap \p {\mathcal M}$,
\begin{align*}
G_n &=\rho e^{\gamma(u)} \sum_{l=1}^n (u_l-\b \varphi_l) \left(u_{ln}
-\beta_n \varphi_l -\b \varphi_{ln}\right)+ G \a_n\\
& =\rho e^{\gamma(u)} \sum_{l=1}^{n-1} u_l u_{nl} + \la G.
\end{align*}
By using again the boundary condition $u_n=\eta(x, u)$ on
$\p\mathcal{M}$ we have
\begin{align*}
u_{nl} &=e_l(u_n)-d u\left(\nabla_{e_l}
e_n\right)=e_l(\eta)-\sum_{k=1}^{n-1} \l\nabla_{e_l} e_n, e_k
\r_g u_k \\
&=\eta_l+\eta_u u_l + \sum_{k=1}^{n-1} \Pi(e_k, e_l) u_k.
\end{align*}
Consequently, by using (\ref{2.1}),  (\ref{1.12}) and the fact
$G(x_0)\ge 1$ we can choose a large universal number $\la$
such that
\begin{align*}
G_n(x_0)&=\rho e^{\gamma(u)} \sum_{k,l=1}^{n-1} \Pi(e_k, e_l) u_k u_l +
\rho e^{\gamma(u)} \sum_{l=1}^{n-1} (\eta_l +\eta_u u_l)u_l +\la G\\
&\ge \left(\la-2C_1-2|\eta_u|\right) G+\rho e^{\gamma(u)} \sum_{l=1}^{n-1} \eta_l u_l\\
&>0.
\end{align*}
However, by the maximality of $G(x_0)$ we have $G_n(x_0)\le 0$. We
thus derive a contradiction.

Therefore $x_0\in \mathcal{O}_1\backslash\p \mathcal{M}$.
Choose normal coordinates around $x_0$ such that
\begin{equation}\label{2.3}
(U_{ij})=\left(u_{ij}+u_i u_j-\frac{1}{2}|\nabla u|_g^2
g_{ij}+(A_g)_{ij}\right)
\end{equation}
is diagonal at $x_0$. Then, by setting $\xi_l=u_l-\beta
\varphi_l$, we have at $x_0$ that
\begin{equation}\label{2.4}
0=G_i=\rho e^{\a} \sum_{l} \left(u_{li}- \b_i
\varphi_l-\b\varphi_{li})\right)\xi_l +\a_i G +\frac{\rho_i}{\rho} G
\end{equation}
and
\begin{align}\label{2.5}
0\ge (G_{ij}) &=\Big( (\a_{ij}-\a_i\a_j) G +\frac{\rho \rho_{ij}
-2\rho_i\rho_j}{\rho^2} G -\frac{\a_i\rho_j+\a_j\rho_i}{\rho} G\nonumber\\
&  + \rho e^{\a} \sum_l\left(u_{lij}-\b_{ij} \varphi_l -\b_i
\varphi_{lj} -\beta_j \varphi_{li} -\beta \varphi_{lij}\right)
\xi_l \nonumber \\
& +\rho  e^{\a} \sum_l\left(u_{li}-\b_i \varphi_l-\beta
\varphi_{li}\right) \left(u_{lj}-\b_j \varphi_l-\beta
\varphi_{lj}\right) \Big).
\end{align}
Let $F^{ij}:=\frac{\p F}{\p U_{ij}}(U)$. It is well-known that
$(F^{ij})$ is positive definite and $F^{ij}=f_i \delta_{ij}$ at $x_0$
(see \cite{CNS}). It then follows from (\ref{2.5}) that
\begin{align}\label{2.6}
0& \ge e^{-\a} F^{ij} G_{ij}\nonumber\\
&=e^{-\a} G\sum_i f_i (\a_{ii}-\a_i^2) + e^{-\a} G\sum_i f_i
\frac{\rho \rho_{ii}-2\rho_i^2}{\rho^2} -2 e^{-\a} G
\sum_i f_i\a_i \frac{\rho_i}{\rho}\nonumber\\
&\quad +\rho\sum_{i,l} f_i \left(
u_{lii}-\b_{ii}\varphi_l-2\b_i
\varphi_{li} -\b \varphi_{lii} \right) \xi_l \nonumber\\
&\quad +\rho\sum_{i,l} f_i \left(u_{li}-\beta_i \varphi_l
-\beta\varphi_{li}\right)^2\nonumber\\
&\ge \frac{1}{2}\rho|\xi|^2 \sum_i f_i (\a_{ii}-\a_i^2) -|\xi|^2
\sum_i f_i \a_i \rho_i -C{\mathcal T} |\xi|^2 +{\mathcal E},
\end{align}
where
$$
{\mathcal E}:=\rho \sum_{i,l} f_i\left(
u_{lii}-\b_{ii}\varphi_l-2\b_i \varphi_{li} -\b \varphi_{lii}
\right) \xi_l.
$$
Since $G(x_0)\ge 1$, we have $|\nabla u|\le C |\xi|$ for some
universal constant $C$.  Note that
$$
\b_i=\eta_i+\eta_u u_i \quad \mbox{and} \quad \b_{ii}
=\eta_{ii}+2\eta_{iu} u_i+ \eta_{uu} u_i^2 +\eta_u u_{ii}.
$$
By using Ricci identity $u_{lii}=u_{iil}+ R_{rili} u_r$,
where $R_{ijkl}$ denotes the Riemann curvature tensor of $g$, we
have
$$
{\mathcal E} \ge \rho \sum_{i,l} f_i u_{iil}\xi_l -\rho \eta_u \sum_{i,l} f_i u_{ii}
\varphi_l \xi_l -C\rho {\mathcal T} |\xi|^3.
$$
By the degree one homogeneity of $f$ and (\ref{2.3}) we have
$$
-\rho \eta_u \sum_{i,l} f_i u_{ii} \varphi_l \xi_l\ge -C\rho
{\mathcal T} |\xi|^3.
$$
Therefore by using (\ref{1.11}),
(\ref{2.4}) and the fact $|\nabla \rho|\le C\sqrt{\rho}$ we have
\begin{align*}
{\mathcal E} &\ge \rho\sum_{i,l} f_i u_{iil} \xi_l -C\rho{\mathcal T}
|\xi|^3\nonumber\\
&=\rho\sum_{i, l} f_i\left(U_{ii}-u_i^2+\frac{1}{2}|\nabla u|^2 g_{ii}
-(A_g)_{ii} \right)_l \xi_l -C\rho{\mathcal T} |\xi|^3\nonumber\\
&\ge -2 \rho\sum_{i, l} f_i u_{il}u_i \xi_l +\rho {\mathcal T} \sum_{k,l}
u_{kl}u_k \xi_l-C\rho {\mathcal T} |\xi|^3\nonumber\\
&\ge -2 \rho \sum_{i, l} f_i \left( u_{li} -\b_i \varphi_l -\b
\varphi_{il}\right) \xi_l u_i \nonumber\\
&\quad  + \rho \mathcal{T}\sum_{k,l}\left(u_{lk}- \b_k \varphi_l -\b
\varphi_{lk}\right)  \xi_l u_k -C\rho {\mathcal T}
|\xi|^3\nonumber\\
&=|\xi|^2 \sum_i f_i u_i(\rho \a_i+\rho_i) -\frac{1}{2} \mathcal{T} |\xi|^2
\sum_k u_k (\rho \a_k +\rho_k) -C\rho\mathcal{T} |\xi|^3\\
&\ge \rho|\xi|^2 \sum_i f_i u_i \a_i -\frac{1}{2}\rho {\mathcal T} |\xi|^2 \sum_k u_k \a_k
-C \sqrt{\rho} {\mathcal T} |\xi|^3.
\end{align*}
Plugging this estimate into (\ref{2.6}) we have
\begin{align*}
0&\ge \frac{1}{2}\rho|\xi|^2 \sum_i f_i (\a_{ii}-\a_i^2)+\rho|\xi|^2 \sum_i
f_i u_i \a_i -|\xi|^2 \sum_i f_i \a_i \rho_i\\
&\quad -\frac{1}{2}\rho |\xi|^2 {\mathcal T}\sum_k u_k\a_k
-C\sqrt{\rho} {\mathcal T}|\xi|^3.
\end{align*}
Note that
$$
\a_i =\la \varphi_i +\gamma'\left((1-\eta_u \varphi)u_i-\eta_i
\varphi-\eta \varphi_i\right)
$$
and
\begin{align*}
\a_{ii}&=\la \varphi_{ii}+\gamma''\left((1-\eta_u
\varphi)u_i-\eta_i \varphi-\eta \varphi_i\right)^2\\
& +\gamma'\left((1-\eta_u \varphi)u_{ii}-\eta_{uu} \varphi u_i^2
-2\eta_{iu} u_i \varphi -2\eta_u u_i \varphi_i -\eta_{ii} \varphi
-2\eta_i \varphi_i -\eta \varphi_{ii}\right).
\end{align*}
It follows that
\begin{align*}
0&\ge \frac{1}{2}\rho |\xi|^2 \sum_i f_i
u_i^2\left\{\left(\gamma''-(\gamma')^2\right) (1-\eta_u \varphi)^2
+2\gamma'(1-\eta_u \varphi)-\gamma'\eta_{uu} \varphi\right\}\\
&\quad +\frac{1}{2}\rho |\xi|^2 \gamma'(1-\eta_u \varphi)\sum_i f_i
u_{ii} -\frac{1}{2}\rho \gamma'{\mathcal T} |\xi|^2 |\nabla u|^2
(1-\eta_u \varphi)-C\sqrt{\rho} {\mathcal T} |\xi|^3
\end{align*}
Note that $u_{ii}=U_{ii}-u_i^2+\frac{1}{2}|\nabla u|^2
g_{ii}-(A_g)_{ii}$. Using the degree one homogeneity of $f$ we
obtain
\begin{align*}
0&\ge \frac{1}{2}\rho |\xi|^2 \sum_i f_i
u_i^2\left\{\left(\gamma''-(\gamma')^2\right) (1-\eta_u \varphi)^2
+\gamma'(1-\eta_u \varphi)-\gamma'\eta_{uu} \varphi\right\}\\
&\quad -\frac{1}{4}\rho \gamma'{\mathcal T} |\xi|^2 |\nabla u|^2
(1-\eta_u \varphi)-C\sqrt{\rho}{\mathcal T} |\xi|^3
\end{align*}
From the definition of $\gamma$ and (\ref{2.2}) one can verify
that
$$
\left(\gamma''-(\gamma')^2\right) (1-\eta_u \varphi)^2
+\gamma'(1-\eta_u \varphi)-\gamma'\eta_{uu} \varphi\ge 0 \quad
\mbox{and} \quad \gamma'\le -c_0<0.
$$
Thus
$$
0\ge c_0\rho {\mathcal T} |\xi|^2|\nabla u|^2-C\sqrt{\rho}
{\mathcal T}|\xi|^3\ge c_0\rho {\mathcal T}|\xi|^4-C\sqrt{\rho}
{\mathcal T}|\xi|^3.
$$
This gives the desired estimate.
\end{proof}

Before giving the proof of Theorem \ref{T2.1}, let us recall the
condition $(H_\a)$ on $(f, \Ga)$ introduced in \cite{LL03}.

\begin{definition}\label{D2.1}
We say $(f, \Ga)$ satisfies condition ($H_\a$) for some $\a>0$ if
there exists some positive constants $\varepsilon_1$ and $c_1$
such that for any $(\la, \xi)\in \Ga\times {\mathbb R}^n$ satisfying
$$
f(\la)\le \a, \quad |\xi|\ge \varepsilon_1^{-1} \quad \mbox{and}
\quad \left|\xi_i\left(\la_i-\frac{1}{2}|\xi|^2\right)\right|\le
\varepsilon_1|\xi|^3 \mbox{ for }1\le i\le n,
$$
there holds
$$
\sum_i
f_i(\la)\left\{\left(\la_i+\frac{1}{2}|\xi|^2-\xi_i^2\right)^2
+\xi_i^2(|\xi|^2-\xi_i^2)\right\}\ge c_1|\xi|^4\sum_i f_i(\la).
$$
\end{definition}

Some discussions have been given in \cite{LL03} on the condition
$(H_\a)$ for $(f, \Ga)$. Here are two remarks.

\begin{Remark}\label{R2.1}
(i) It is easy to check that if $f$ is homogeneous on $\Ga$, then
$(f, \Ga)$ satisfies the condition $(H_1)$ if and only if $(f,
\Ga)$ satisfies the condition $(H_\a)$ for each $\a>0$.

(ii) From the proof of \cite[Lemma 2.4]{GW03} and \cite[Theorem
3]{GW04}, one can see that the following two classes of $(f, \Ga)$
satisfy the condition $(H_\a)$:
$$
(f, \Ga)=(\sigma_k^{1/k}, \Ga_k) \quad \mbox{with } 1\le k\le n
$$
and
$$
(f, \Ga)=\left((\sigma_k/\sigma_l)^{1/(k-l)}, \Ga_k\right)
$$
with $0\le l<k\le n$ and $(n-k+1)(n-l+1)>2(n+1)$.
\end{Remark}

\begin{proof}[Proof of Theorem \ref{T2.1}]  Consider the function
$$
G:=\frac{1}{2}\rho \mbox{e}^{\la \varphi} |\nabla u-c e^{-u}\nabla
\varphi|_g^2\quad \mbox{on } {\mathcal O}_1,
$$
where $\varphi$ and $\rho$ are as in the proof of Theorem \ref{T1.3},
and $\la$ is a fixed positive constant satisfying
$$
\la>2 C_1+ 2 |c| e^{C_0}.
$$
We need to show that $G$ can be bounded by some universal constant
$C$ depending only on $n$, $c$, $(f, \Ga)$, $g$, $C_0$, ${\mathcal O}_1$ and
${\mathcal O}_2$. Suppose the maximum of $G$ over
$\mathcal{O}_1$ is attained at some point $x_0\in \mathcal{O}_1$. In
the following we will always assume that $G(x_0)\ge 1$.

Similar to the proof of Theorem \ref{T1.3} we have $x_0\in
\mathcal{O}_1\backslash\p \mathcal{M}$. As before we choose normal
coordinates around $x_0$ such that $(U_{ij})$ is diagonal at
$x_0$. Then, by setting $\xi_l=u_l-c e^{-u} \varphi_l$, we have at
$x_0$ that
\begin{equation}\label{2.8}
0=G_i=\rho e^{\la\varphi} \sum_{l} \left(u_{li}-
ce^{-u}(\varphi_{li}-\varphi_l u_i)\right)\xi_l +\la G \varphi_i +
\frac{\rho_i}{\rho} G
\end{equation}
and
\begin{align}\label{2.9}
0\ge (G_{ij}) &=\Big( (\la \varphi_{ij}-\la^2\varphi_i
\varphi_j) G +\frac{\rho\rho_{ij}-2\rho_i\rho_j}{\rho^2} G
-\frac{\rho_i \varphi_j+\rho_j\varphi_i}{\rho} \la G\nonumber\\
&  + \rho e^{\la\varphi} \sum_l\left(u_{lij}-ce^{-u}(\varphi_{lij}
-\varphi_{li}u_j-\varphi_{lj}u_i-\varphi_l u_{ij} +\varphi_l u_i
u_j)\right) \xi_l \nonumber \\
& + \rho e^{\la\varphi}
\sum_l\left(u_{li}-ce^{-u}(\varphi_{li}-\varphi_l u_i)\right)
\left(u_{lj}-ce^{-u}(\varphi_{lj}-\varphi_l u_j)\right) \Big).
\end{align}
It then follows from (\ref{2.9}) that
\begin{equation}\label{2.10}
0\ge e^{-\la\varphi} F^{ij} G_{ij} =I
+II -\frac{C}{\rho}  \mathcal{T} G.
\end{equation}
where
\begin{align*}
I &:=\rho\sum_{i, l}f_i \left(u_{lii}-
ce^{-u}(\varphi_{lii}-2\varphi_{li}u_i-\varphi_l u_{ii} +\varphi_l
u_i^2)\right)\xi_l,\\
II&:= \rho\sum_{i, l} f_i
\left(u_{li}-ce^{-u}(\varphi_{li}-\varphi_l
u_i)\right)^2,
\end{align*}
Similar to the estimate for ${\mathcal E}$ in the proof of Theorem
\ref{T1.3}, we have by (\ref{2.8}) that
\begin{align}\label{2.11}
I \ge \rho \sum_{i,l}f_i u_{iil} \xi_l -\frac{C}{\sqrt{\rho}} {\mathcal T}
G^{\frac{3}{2}}\ge -\frac{C}{\sqrt{\rho}} {\mathcal T} G^{\frac{3}{2}}.
\end{align}
For the term $II$, by using the elementary inequality
\begin{equation}\label{2.12}
(a+b)^2\ge \frac{1}{2}a^2-b^2\quad \mbox{for any } a, b\in {\mathbb
R}
\end{equation}
we have
\begin{align*}
II &\ge \frac{1}{2}\rho \sum_{i,l} f_i (u_{li}+(A_g)_{li})^2 -
\rho \sum_{i,l} f_i \left((A_g)_{il}+c
e^{-u}(\varphi_{li}-\varphi_l u_i)\right)^2 \nonumber\\
&\ge \frac{1}{2} \rho \sum_{i,l} f_i (u_{li}+(A_g)_{li})^2- C{\mathcal
T} G.
\end{align*}
Using (\ref{2.3}) we then obtain
\begin{align}\label{2.13}
II &=\frac{1}{2}\rho \sum_{i,l}
f_i\left(U_{li}+\frac{1}{2}|\nabla
u|_g^2 g_{li}-u_l u_i\right)^2 -C{\mathcal T} G\nonumber\\
&=\frac{1}{2} \rho \sum_i f_i \left\{\left(U_{ii} +\frac{1}{2}|\nabla
u|_g^2-u_i^2\right)^2+u_i^2\left(|\nabla
u|_g^2-u_i^2\right)\right\} -C{\mathcal T} G.
\end{align}
Let $\xi=(\xi_1, \cdots, \xi_n)$. By using the elementary
inequality (\ref{2.12}) once again, we have
\begin{align*}
\sum_i f_i & \left\{\left(U_{ii}+\frac{1}{2}|\nabla
u|_g^2-u_i^2\right)^2+u_i^2(|\nabla u|_g^2-u_i^2)\right\}\\
&\qquad \ge \frac{1}{2}\sum_i f_i
\left\{\left(U_{ii}+\frac{1}{2}|\xi|^2-\xi_i^2\right)^2
+\xi_i^2(|\xi|^2-\xi_i^2)\right\}\\
&\qquad\quad -\sum_i f_i \left(\frac{1}{2}|\nabla
u|_g^2-\frac{1}{2}|\xi|^2 -u_i^2+\xi_i^2\right)^2 \\
&\qquad \quad+\sum_i f_i \left(u_i^2(|\nabla
u|_g^2-u_i^2)-\xi_i^2(|\xi|^2-\xi_i^2)\right)
\end{align*}
It is easy to see that
$$
\left|\frac{1}{2}|\nabla u|_g^2 -\frac{1}{2}|\xi |^2 -u_i^2+
\xi_i^2\right|\le C(1+|\nabla u|_g)\le C|\xi|
$$
and
$$
\left|u_i^2(|\nabla
u|_g^2-u_i^2)-\xi_i^2(|\xi|^2-\xi_i^2)\right|\le C(1 +|\nabla
u|_g^3)\le C|\xi|^3.
$$
Consequently
\begin{align*}
 \sum_i  f_i&  \left\{\left(U_{ii}+\frac{1}{2}|\nabla
u|_g^2-u_i^2\right)^2+u_i^2(|\nabla u|_g^2-u_i^2)\right\}\\
&\ge \frac{1}{2} \sum_i f_i
\left\{\left(U_{ii}+\frac{1}{2}|\xi|^2-\xi_i^2\right)^2
+\xi_i^2(|\xi|^2-\xi_i^2)\right\}-C{\mathcal T} |\xi|^3.
\end{align*}
This together with (\ref{2.13}) implies that
\begin{equation}\label{2.14}
II\ge \frac{1}{4} \rho \sum_i f_i
\left\{\left(U_{ii}+\frac{1}{2}|\xi|^2
-\xi_i^2\right)^2+\xi_i^2(|\xi|^2-\xi_i^2)\right\}-\frac{C}{\sqrt{\rho}}{\mathcal T}
G^{\frac{3}{2}}.
\end{equation}
Combining (\ref{2.10}), (\ref{2.11}) and (\ref{2.14}) yields
\begin{equation}\label{2.15}
\rho \sum_i f_i \left\{\left(U_{ii}+\frac{1}{2}|\xi|^2
-\xi_i^2\right)^2+\xi_i^2(|\xi|^2-\xi_i^2)\right\}\le \frac{C}{\sqrt{\rho}}{\mathcal
T} G^{\frac{3}{2}}+\frac{C}{\rho} {\mathcal T} G.
\end{equation}
In order to apply the ($H_\a$) condition on $(f, \Ga)$ with
$\la_i=U_{ii}$ and $\xi_i=u_i-c e^{-u} \varphi_i$, we need to
check
\begin{equation}\label{2.16}
\left|\xi_i\left(U_{ii}-\frac{1}{2}|\xi|^2\right)\right|\le
\varepsilon_1 |\xi|^3.
\end{equation}
To see this, recall that $(U_{ij})$ is
diagonal, one has
\begin{align*}
\xi_i\left(U_{ii}-\frac{1}{2}|\xi|^2\right) &=\sum_l \xi_l
U_{il}-\frac{1}{2} \xi_i |\xi|^2\\
&  =\sum_l \xi_l \left(u_{il}+u_iu_l-\frac{1}{2}|\nabla u|^2
\delta_{il}+(A_g)_{il}\right)-\frac{1}{2}\xi_i |\xi|^2
\end{align*}
By using (\ref{2.8}) we have
$$
\left|\sum_l \xi_l u_{il}\right|\le \left|\sum_l
\xi_l\left(u_{li}-c e^{-u}(\varphi_{li}-\varphi_l
u_i)\right)\right| +C |\xi|^2\le \frac{C}{\sqrt{\rho}}|\xi|^2.
$$
Moreover, by direct calculation one can see that
$$
\left|\sum_l \xi_l\left(u_i u_l-\frac{1}{2}|\nabla
u|_g^2\delta_{il}\right) -\frac{1}{2}\xi_i|\xi|^2\right| \le
C(1+|\nabla u|_g^2)\le C |\xi|^2.
$$
Consequently we have
$$
\left|\xi_i\left(U_{ii}-\frac{1}{2}|\xi|^2\right)\right|\le
\frac{C}{\sqrt{\rho}} |\xi|^2\le \frac{C}{\sqrt{G}} |\xi|^3 \le
\varepsilon_1 |\xi|^3.
$$
if we further assume that $G(x_0)\ge C^2/\varepsilon_1^2$.
Therefore we may apply the ($H_\a$) condition on $(f, \Ga)$ to
(\ref{2.15}) to get
$$
c_1\rho\mathcal{T} |\xi|^4 \le \frac{C}{\sqrt{\rho}} {\mathcal T}
G^{\frac{3}{2}}+\frac{C}{\rho}{\mathcal T} G\le \frac{C}{\rho} {\mathcal T}
\left(G^{\frac{3}{2}}+G\right).
$$
Consequently $G^2\le C(G^{\frac{3}{2}}+G)$
which implies $G(x_0)\le C$ for some universal constant $C$.
\end{proof}

\section{\bf $C^2$ estimates: general equations}
\setcounter{equation}{0}

The aim of this section is to show Theorem \ref{T1.5}. We may assume 
$\mathcal{O}_1\cap \p {\mathcal M}\ne
\emptyset$ since otherwise the results follow from the well-known
local interior estimates (see \cite{GW03, LL03}). Without loss of
generality, we may also assume that $g$ is conformally flat on
$\mathcal{O}_1$, i.e. there exists a function $\varphi\in
C^\infty({\mathcal O}_1)$ such that $e^{2\varphi} g$ is a flat
metric on ${\mathcal O}_1$. Since ${\mathcal O}_1\cap \p {\mathcal
M}$ is umbilic in $g$, it is also umbilic in the flat metric.
Therefore $\mathcal{O}_1\cap \p {\mathcal M}$ is either a part of
a hyperplane or a part of a sphere in ${\mathbb R}^n$. By
shrinking ${\mathcal O}_1$ and using the conformal diffeomorphism
in ${\mathbb R}^n$ if necessary we may assume that
$$
{\mathcal O}_1=B_4^+:=\{x=(x_1, \cdots, x_n)\in {\mathbb R}^n : |x|<4
\mbox{ and } x_n\ge 0 \}
$$
and
$$
\mathcal{O}_1\cap \p {\mathcal M}:=\p B_4^+\cap \{x_n=0\}.
$$
Observe that $e^{2(u+\varphi)}\tilde{g}$ is also a flat metric, By
using the conformal invariance the function $v:=u+\varphi$
satisfies the following problem
\begin{equation}\label{3.2}
\left\{\begin{array}{lll} F(V):=f(\la(V))=\tilde{\psi}(x, v),\quad
\la(V)\in \Ga \quad \mbox{in } B_4^+,\\
\\
v_n=\tilde{\eta}(x, v) \quad \mbox{ on } \p B_4^+\cap \{x_n=0\},
\end{array}\right.
\end{equation}
where $\tilde{\psi}(x, v)=\psi(x, v-\varphi) e^{-2\varphi}$,
$\tilde{\eta}(x, v)=\eta(x, v-\varphi)e^{-2\varphi}$, $V$ denotes
the matrix function
$$
V:=D^2 v+dv\otimes dv-\frac{1}{2}|Dv|^2 I,
$$
$I$ is the $n\times n$ identity matrix, $D$ is the standard
connection on ${\mathbb R}^n$, $D v$ and $D^2 v$ denote the
gradient and Hessian of $v$ respectively, and $\la(V)$ denote the
eigenvalues of $V$.

Recall that we can assume $h_g=0$ on $\p {\mathcal M}$. So
$\varphi_n=0$ on $\p B_4^+\cap \{x_n=0\}$. Therefore when
$\eta\equiv 0$ and $\frac{\p \psi}{\p \nu}\equiv 0$ on $\p
{\mathcal M}\times {\Bbb R}$, we have $\tilde{\eta}\equiv 0$ and
$\tilde{\psi}_n\equiv 0$ on $\p B_4^+\cap \{x_n=0\}$. Thus, in
order to show Theorem \ref{T1.5}, it suffices to prove the
following result.

\begin{Theorem}\label{T3.2}
Assume that $(f, \Ga)$ satisfies (\ref{1.1})--(\ref{1.6}) and that
$\tilde{\psi}$, $\tilde{\eta}$ satisfy one of the following
conditions:
\begin{enumerate}
\item[(i)] $\tilde{\eta}\equiv 0$ and $\tilde{\psi}_n\equiv 0$ on
$(\p B_4^+\cap \{x_n=0\})\times {\Bbb R}$ and $\tilde{\psi}\in
C^3(B_4^+\times{\Bbb R})$;

\item[(ii)] $\tilde{\eta}$ is positive and $\tilde{\psi}$ is any
function.
\end{enumerate}
If $v\in C^4(B_4^+)$ is a solution of (\ref{3.2}) satisfying
\begin{equation}\label{3.3}
|v|\le C_0 \quad \mbox{and} \quad  |Dv|\le C_0 \quad \mbox{in }
B_4^+
\end{equation}
for some positive constant $C_0$, then
\begin{equation}\label{3.4}
|D^2 v|\le C\quad  \mbox{in } B_1^+
\end{equation}
for some positive constant $C$ depending only on $n$, $C_0$, $(f,
\Ga)$, $\tilde{\psi}$ and $\tilde{\eta}$.
\end{Theorem}

\begin{proof}[Proof of Theorem \ref{T3.2} under condition (i)]
This case can be reduced to the local interior estimates. To this
end, we define
\begin{equation}\label{3.5}
\bar{v}(x', x_n)=\left\{\begin{array}{lll} v(x', x_n), &
\mbox{when } x_n\ge 0,\\
\\
v(x', -x_n), & \mbox{when } x_n\le 0
\end{array}\right.
\end{equation}
and
$$
\bar{\psi}((x',x_n), t)= \left\{\begin{array}{lll}
\tilde{\psi}((x',x_n), t), & \mbox{when } x_n\ge 0, ~t\in {\Bbb
R},\\
\\
\tilde{\psi}((x', -x_n), t), & \mbox{when } x_n\le 0, ~ t\in {\Bbb
R}.
\end{array}\right.
$$
Since $v_n(x', 0)=0$ and $\tilde{\psi}_n((x',0),t)=0$, it is easy
to see $\tilde{v}\in C^2(B_4)$ and $\bar{\psi}\in
C^{2,1}(B_4\times {\Bbb R})$. Let
$$
\overline{V}:=D^2\bar{v}+d\bar{v}\otimes d\bar{v}-\frac{1}{2}|D
\bar{v}|^2 I.
$$
By direct calculation one can see that
$$
\overline{V}(x', x_n)=\left\{\begin{array}{lll} V(x', x_n),&
\mbox{when } x_n\ge 0,\\
\\
{\mathcal Q}^t V(x', -x_n) \mathcal{Q}, &\mbox{when } x_n\le 0,
\end{array}\right.
$$
where ${\mathcal Q}$ is the orthogonal matrix ${\mathcal
Q}:=\mbox{diag}[1, \cdots, 1, -1]$. Therefore
\begin{equation}\label{3.6}
F(\overline{V}):=f(\la(\overline{V}))=\bar{\psi}(x, \bar{v}),
\quad \la(\overline{V})\in \Ga \quad \mbox{in } B_4.
\end{equation}
It then follows from \cite[Lemma 17.16]{GT01} that $\bar{v}\in
C^{4, \a}(B_2)$ for any $\a\in (0,1)$. Now the local interior
estimates in \cite{LL03} can be applied to obtain $|D^2
\bar{v}|\le C$ in $B_1$ for some constant $C$ depending only on
$n$, $C_0$, $(f, \Ga)$ and $\tilde{\psi}$. This in particular
implies (\ref{3.4}).
\end{proof}

\begin{Remark}\label{R3.1}
Note that the function $\bar{v}$ defined by (\ref{3.5}) satisfies
(\ref{3.6}) and has uniform $C^2$ estimates. Since $f$ is concave,
by Evans-Krylov theory and Schauder theory we can obtain uniform
estimates on $\|\bar{v}\|_{C^{4, \a}(B_{1/2})}$ for any $\a\in (0,
1)$. Therefore, for a solution $u\in C^2({\mathcal M})$ of
(\ref{1.11}) with $\psi$ and $\eta$ satisfying (i) in Theorem
\ref{T1.5}, if (\ref{3.1}) holds then
$$
\|u\|_{C^{4, \a}({\mathcal O}_2)}\le C
$$
for some constant $C$ depending only on $n$, $C_0$, $\a$, $(f,
\Ga)$, $\psi$,  ${\mathcal O}_1$ and ${\mathcal O}_2$.
\end{Remark}

Next we will prove Theorem \ref{T3.2} under condition (ii).  We
will use $\{e_1, \cdots, e_n\}$ to denote the standard orthonormal
basis in ${\mathbb R}^n$, i.e. $e_i=(0, \cdots, 1, \cdots, 0)$,
where $1$ is in the $i$th spot and $0$ elsewhere. The following
result gives the double tangential derivative estimates without
any restrictions on $\tilde{\psi}$ and $\tilde{\eta}$.

\begin{Lemma}\label{L3.1}
Assume that $(f, \Ga)$ satisfies (\ref{1.1})--(\ref{1.6}). If
$v\in C^4(B_4^+)$ is a solution of (\ref{3.2}) satisfying
\begin{equation}\label{3.7}
|v|\le C_0 \quad \mbox{and} \quad  |Dv|\le C_0 \quad \mbox{in }
B_4^+
\end{equation}
for some positive constant $C_0$, then there exists a positive
constant $C_1$ depending only on $n$, $C_0$, $(f, \Ga)$,
$\tilde{\psi}$ and $\tilde{\eta}$ such that
\begin{equation}\label{3.8}
v_{\tau\tau}\le C_1\quad  \mbox{in } B_3^+
\end{equation}
for any vector $\tau\in \mbox{span}\{e_1, \cdots, e_{n-1}\}$ with
$|\tau|=1$.
\end{Lemma}
\begin{proof}
By rotation it suffices to establish (\ref{3.8}) for $\tau=e_1$.
Let $\rho\in C^\infty_0(B_4)$ be a radial cut-off function such
that $0\le \rho\le 1$ in ${\mathbb R}^n$, $\rho=1$ in $B_3$, and
$|D\rho|\le C\rho^{\frac{1}{2}}$  in $B_4$. Consider the function
$$
H=\rho e^{\beta x_n} \left( v_{11}+ v_1^2\right) \quad \mbox{on }
B_4^+,
$$
where $\beta$ is a large positive constant to be determined later.
Suppose the maximum of $H$ over $B_4^+$ is attained at some point
$x_0$, then either $x_0\in B_4^+\cap\{x_n>0\}$ or $x_0 \in B_4\cap
\{x_n=0\}$. In the following we will always assume that $H(x_0)\ge
1$ and $v_{11}(x_0)\ge 1$; otherwise we are done.

Let us first calculate $H_n$ on $B_4\cap \{x_n=0\}$. Since $\rho$
is radially symmetric, we have $\rho_n=0$ on $x_n=0$. Thus, by
using the boundary condition $v_n=\tilde{\eta}(x, v)$,  it is easy
to see that
\begin{align*}
H_n&=\rho\left(v_{11n}+2v_1v_{1n}+ \beta(v_{11}+ v_1^2)\right)\\
&=\rho\left((\beta+\tilde{\eta}_v)v_{11}+ (\beta+\tilde{
\eta}_{vv}+2\tilde{\eta}_v) v_1^2+ (2\tilde{\eta}_1+
2\tilde{\eta}_{1v}) v_1 +\tilde{\eta}_{11}\right).
\end{align*}
If $x_0\in B_4\cap \{x_n=0\}$, then, using (\ref{3.7}) and
$v_{11}(x_0)\ge 1$, we have $H_n(x_0)>0$ by choosing $\beta$ large
enough. But by the maximality of $H(x_0)$ we have $H_n(x_0)\le 0$.
We thus derive a contradiction. Therefore $x_0\in
B_4^+\cap\{x_n>0\}$.

Now at $x_0$ we have
\begin{equation}\label{3.9}
0=H_i=\left(\frac{\rho_i}{\rho} +\beta \delta_{in}\right) H +\rho
e^{\beta x_n}\left(v_{11i}+2v_1 v_{1i}\right), \quad 1\le i\le n
\end{equation}
and
\begin{align*}
0\ge (H_{ij})&=\Big(\frac{\rho \rho_{ij}-2\rho_i\rho_j}{\rho^2} H
- \beta^2 \delta_{in}\delta_{jn} H-\frac{\rho_i \delta_{jn}+\rho_j
\delta_{in}}{\rho} \beta H \\
&\quad + \rho e^{\beta x_n} (v_{11ij}+2v_1 v_{1ij} + 2
v_{1i}v_{1j})\Big ).
\end{align*}
Let $F^{ij}:=\frac{\p F}{\p V_{ij}}(V)$ and ${\mathcal T}:=\sum_i
F^{ii}$. We know that $(F^{ij})$ is positive definite and
${\mathcal T}\ge f(1, \cdots, 1)>0$, Thus we have at $x_0$ that
\begin{align}\label{3.10}
0&\ge e^{-\beta x_n}  F^{ij} H_{ij}\nonumber\\
&=e^{-\beta x_n} H F^{ij}\frac{\rho
\rho_{ij}-2\rho_i\rho_j}{\rho^2} -\beta^2 e^{-\beta x_n} H F^{nn}
-2\beta e^{-\beta x_n} H F^{in} \frac{\rho_i}{\rho}\nonumber\\
&\quad + \rho F^{ij} \left(v_{11ij} +2v_1
v_{1ij}+2v_{1i}v_{1j}\right)\nonumber\\
&\ge -\frac{C}{\rho}{\mathcal T} H+\rho F^{ij} \left(v_{11ij}+2v_1
v_{1ij}+2v_{1i}v_{1j}\right)
\end{align}
By differentiating the equation (\ref{3.2}) twice and using the
concavity of $F$ we have
$$
F^{ij} \Big(v_{ijk}+2v_{ik} v_j- \sum_l v_l v_{lk}
\delta_{ij}\Big)=\tilde{\psi}_k +\tilde{\psi}_v v_k, \quad 1\le
k\le n
$$
and
\begin{align*}
F^{ij}\Big( v_{ij11}+2v_{i11} v_j +2 v_{i1} v_{j1} &-\sum_l (v_l
v_{l11}+ v_{l1}^2) \delta_{ij}\Big)  \\
&\ge \tilde{\psi}_{11}+2\tilde{\psi}_{1v} v_1+\tilde{\psi}_{vv}
v_1^2 +\tilde{\psi}_v v_{11}.
\end{align*}
Therefore, by using (\ref{3.9}) and (\ref{3.7}),
\begin{align}\label{3.11}
\rho F^{ij} v_{ij11}&\ge -2\rho F^{ij}\left(v_{i11}v_j+ v_{i1}
v_{j1}\right) +\rho {\mathcal T} \sum_l \left(v_l
v_{l11}+v_{l1}^2\right)-\frac{C}{\rho} {\mathcal T} H\nonumber\\
&\ge 4\rho v_1 F^{ij} v_{1i} v_j -2\rho {\mathcal T} v_1 \sum_l
v_l v_{l1} -2\rho F^{ij} v_{i1}v_{j1}\nonumber\\
&\quad +\rho {\mathcal T} \sum_l v_{l1}^2-\frac{C}{\rho} {\mathcal
T} H.
\end{align}
By using (\ref{3.2}) and (\ref{3.7}) we have
\begin{align}\label{3.12}
2\rho  F^{ij} v_1 v_{1ij}&=2\rho v_1  F^{ij}
\Big(V_{ij}-v_iv_j+\frac{1}{2}|Dv|^2\delta_{ij}\Big)_1\nonumber\\
&=2\rho v_1(\tilde{\psi}_1+\tilde{\psi}_v v_1)-4\rho v_1  F^{ij}
v_{i1}v_j
+2\rho v_1 {\mathcal T} \sum_l v_{l1} v_l\nonumber\\
&\ge -4\rho v_1  F^{ij} v_{i1}v_j +2\rho v_1{\mathcal T} \sum_l
v_{l1} v_l-\frac{C}{\rho}{\mathcal T} H.
\end{align}
Combining (\ref{3.10})--(\ref{3.12})  yields
$$
0\ge -\frac{C}{\rho}{\mathcal T} H+\rho {\mathcal T} \sum_l
v_{l1}^2 \ge -\frac{C}{\rho}{\mathcal T} H+\rho{\mathcal T}
v_{11}^2 \ge -\frac{C}{\rho}{\mathcal T} H+\frac{1}{\rho}
{\mathcal T} H^2.
$$
Consequently, we have $H(x_0)\le C$, and the proof is complete.
\end{proof}

\begin{Lemma}\label{L3.2}
Under the hypotheses of Theorem \ref{T3.2} with (ii) satisfied,
there exists a positive constant $C$ depending only on $n$, $C_0$,
$(f, \Ga)$, $\tilde{\psi}$ and $\tilde{\eta}$ such that
\begin{equation}\label{3.15}
|v_{nn}(x',0)|\le C\quad  \mbox{ whenever } |x'|\le 2.
\end{equation}
\end{Lemma}
\begin{proof}
It is convenient to consider the function $w:=e^v$. By direct
calculation and the degree one homogeneity of $F$ it follows from
(\ref{3.2}) that $w$ satisfies
\begin{equation}\label{3.16}
\left\{\begin{array}{lll} F(W):=f(\la(W)) =\hat{\psi}(x, w),
\quad \la(W)\in \Ga \quad \mbox{in } B_4^+\\
 \\
w_n=\hat{\eta}(x, w) \quad \mbox{ on } \p B_4^+\cap \{x_n=0\}.
\end{array}\right.
\end{equation}
where $\hat{\psi}(x,w)=w\tilde{\psi}(x, \log w)$, $\hat{\eta}(x,
w)=w\tilde{\eta}(x, \log w)$ and
$$
W:=D^2 w -\frac{1}{2w} |D w|^2 I.
$$
Moreover, by using (\ref{3.3}) one can see that there is a
positive constant $C_2$ depending only on $C_0$ such that
\begin{equation}\label{3.17}
\frac{1}{C_2}\le w\le C_2 \quad \mbox{and} \quad |D w|\le C_2
\quad \mbox{in } B_4^+.
\end{equation}

We now introduce a linear elliptic differential operator
${\mathcal L}$ on $B_4^+$ by
$$
{\mathcal L}\phi = F^{ij}\phi_{ij} - w^{-1} {\mathcal T} \sum_l
w_l \phi_l, \quad \forall \phi\in C^2(B_4^+).
$$
where $F^{ij}:=\frac{\p F}{\p W_{ij}}(W)$ and ${\mathcal
T}:=\sum_i F^{ii}$.  Recall that $(F^{ij})$ is positive definite
and ${\mathcal T}\ge f(1, \cdots, 1)>0$. By differentiating the
equation (\ref{3.16}) with respect to $x_n$ we obtain
$$
F^{ij}\left( w_{ijn}-w^{-1}w_l w_{ln}\delta_{ij}+\frac{1}{2}w^{-2}
w_n |D w|^2 \delta_{ij}\right)=\hat{\psi}_n+\hat{\psi}_w w_n.
$$
Thus
$$
{\mathcal L} w_n=\hat{\psi}_n +\hat{\psi}_w w_n -\frac{1}{2}
w^{-2} w_n {\mathcal T} |D w|^2.
$$
By using (\ref{3.17}) and the degree one homogeneity of $f$ we
have
$$
|{\mathcal L}\hat{\eta}|\le |\hat{\eta}_w F^{ij}
w_{ij}|+C{\mathcal T}\le C{\mathcal T} \quad \mbox{in } B_4^+.
$$
Consequently
\begin{equation}\label{3.18}
\left|{\mathcal L}\left(w_n -\hat{\eta}(x, w)\right)\right| \le
C{\mathcal T} \quad \mbox{in } B_4^+.
\end{equation}

Next we consider the function $w_n$. Since $\la(W)\in \Ga\subset
\Ga_1$, we have $\Delta w\ge 0$ in $B_4^+$. Therefore from Lemma
\ref{L3.1} it follows that
$$
w_{nn}\ge -\sum_{i=1}^{n-1} w_{ii} \ge -(n-1)C_1 \quad \mbox{in }
B_3^+.
$$
Since $\hat{\eta}$ is positive, we have $\hat{\eta}\ge 2\a_0$ on
$\p B_4^+\cap \{x_n=0\}$ for some universal constant $\a_0>0$.
Therefore
$$
w_n(x', x_n)+ (n-1)C_1 x_n\ge w_n(x', 0)\ge 2\a_0 \quad \mbox{for
} (x', x_n) \in B_3^+.
$$
Thus there exists a universal constant $0<\varepsilon_0\le 1$ such
that
\begin{equation}\label{3.19}
w_n(x', x_n)\ge \a_0>0 \quad \mbox{for } (x', x_n)\in B_3^+\cap
\{x_n\le \varepsilon_0\}.
\end{equation}

In order to establish (\ref{3.15}), it suffices to show that
$|w_{nn}(0)|\le C$ for some universal constant $C$. Consider the
function
$$
\phi=A x_n +B|x|^2\pm \left(w_n-\hat{\eta}(x, w)\right) \quad
\mbox{on } \overline{B_{\varepsilon_0}^+},
$$
where $A$ and $B$ are sufficiently large positive constants to be
chosen below. Clearly on $\p B_{\varepsilon_0}^+\cap \{x_n=0\}$ we
have $\phi\ge 0$ since $w_n-\hat{\eta}(x, w)=0$ there. Also, by
using (\ref{3.17}) we may choose $B$ large enough so that $\phi\ge
0$ on $\p B_{\varepsilon_0}^+\cap \{x_n>0\}$. Thus
\begin{equation}\label{3.20}
\phi\ge 0 \quad \mbox{on } \p B_{\varepsilon_0}^+.
\end{equation}
With the number $B$ chosen above, by using (\ref{3.17}),
(\ref{3.18}) and (\ref{3.19}) we have
\begin{align}\label{3.21}
{\mathcal L} \phi&= -A w^{-1} w_n {\mathcal T}+ B{\mathcal
L}(|x|^2)\pm {\mathcal L}\left(w_n-\hat{\eta}(x, w)\right)\nonumber\\
&\le -A C_1^{-1}\a_0 {\mathcal T}+ C {\mathcal T}\nonumber\\
&<0
\end{align}
in $B_{\varepsilon_0}^+$ if we choose $A$ large enough.

By the maximum principle, it follows from (\ref{3.20}) and
(\ref{3.21}) that $\phi\ge 0$ in $B_{\varepsilon_0}^+$. Since
$\phi(0)=0$, we therefore have $\phi_n(0)\ge 0$. consequently
$|w_{nn}(0)|\le C$ for some universal constant $C$.
\end{proof}

\begin{proof}[Proof of Theorem \ref{T3.2} under condition (ii)]
Consider the function
$$
\widetilde{G}(x, \xi)=\rho(x) \left(v_{\xi \xi}(x)+\l D v(x),
\xi\r^2\right), \qquad \forall x\in \overline{B_2^+} \mbox { and }
\xi\in {\mathbb S}^n,
$$
where $\rho\in C^\infty_0(B_2)$ is a radial cut-off function such
that $0\le \rho\le 1$ in ${\mathbb R}^n$, $\rho=1$ on $B_1$, and
$|D\rho|\le C\rho^{\frac{1}{2}}$ in $B_2$. Suppose the maximum of
$\widetilde{G}$ over $\overline{B_2^+}\times {\mathbb S}^n$ is
attained at $(\bar{x}, \bar{\xi})$, then either $\bar{x}\in
B_2^+\cap \{x_n>0\}$ or $\bar{x}\in B_2\cap \{x_n=0\}$.

If $\bar{x}\in B_2^+\cap\{x_n>0\}$, then similar to the proof of
\cite[Theorem 1.20]{LL03} one can show that
$\widetilde{G}(\bar{x}, \bar{\xi}) \le C$ for some universal
constant $C$. If $\bar{x}\in B_2\cap \{x_n=0\}$, we may write
$\bar{\xi}=\a e_n +\b \tau$, where $\tau$ is a unit vector in
$\mbox{span}\{e_1, \cdots, e_{n-1}\}$, and $\a$ and $\b$ are two
numbers satisfying $\a^2+\b^2=1$. Then
$$
v_{\bar{\xi} \bar{\xi}}=\a^2 v_{nn} +\b^2 v_{\tau \tau} +2\a \b
v_{n \tau}.
$$
Since $v_n=\tilde{\eta}(x, v)$ on $x_n=0$, we have $v_{n
\tau}=\tilde{\eta}_\tau+\tilde{\eta}_v v_{\tau}$ which is bounded.
Therefore it follows from Lemma \ref{L3.1} and Lemma \ref{L3.2}
that $v_{\bar{\xi}\bar{\xi}}\le C$ at $\bar{x}$ and hence
$\widetilde{G}(\bar{x}, \bar{\xi})\le C$ for some universal
constant $C$.

The above argument shows that $v_{\xi\xi}\le C$ in $B_1^+$ for any
unit vector $\xi\in {\mathbb S}^n$. Since $\Delta v\ge 0$, we also
have $v_{\xi\xi}\ge -C$ in $B_1^+$ for any $\xi\in {\mathbb S}^n$.
Therefore $|D^2v|\le C$ in $B_1^+$.
\end{proof}

\section{\bf $C^2$ estimates: Monge-Ampere type equations}
\setcounter{equation}{0}

In this section we prove Theorem \ref{T4.1}. We first have the following
double normal derivative estimates.

\begin{Lemma}\label{L4.1}
Under the hypotheses of Theorem \ref{T4.1}, there holds
$$
|\nabla^2 u(\nu, \nu)|\le C \quad \mbox{on } \p {\mathcal O}_2\cap
{\mathcal M}
$$
for some positive constant $C$ depending only on $n$, $c$, $g$,
$C_0$, ${\mathcal O}_1$ and ${\mathcal O}_2$.
\end{Lemma}
\begin{proof} Recall that we assume $h_g=0$ on $\p {\mathcal
M}$ and that $d_g(x, \p\mathcal{M})$ is smooth in ${\mathcal
M}_{\delta_0}:=\{x\in {\mathcal M}: d_g(x, \p \mathcal{M})\le
\delta_0\}$.  We may extend the unit inward normal vector field
$\nu$ to a smooth vector field in ${\mathcal M}_{\delta_0}$, still
denoted it as $\nu$, by parallel translating along the unit-speed
geodesics perpendicular to $\p {\mathcal M}$. Clearly $\nabla_\nu
\nu=0$ in ${\mathcal M}_{\delta_0}$. Thus along any such geodesic
$\gamma$ starting from a point $\gamma(0)\in {\mathcal O}_1\cap \p
{\mathcal M}$ we have
$$
\frac{d}{dt} \left(u_{\nu}(\gamma(t))\right)=
\nu(u_{\nu})=\nabla^2 u(\nu, \nu)+ d u(\nabla_\nu \nu)= \nabla^2
u(\nu, \nu).
$$
Since $\la_g(U)\in \Gamma_n$, we have
$$
\nabla^2 u(\nu, \nu)+u_{\nu}^2-\frac{1}{2}|\nabla u|_g^2+A_g(\nu,
\nu)>0.
$$
Therefore it follows from (\ref{4.2}) that there is a universal
constant $C_1$ such that
$\frac{d}{dt}\left(u_\nu(\gamma(t))\right) \ge -C_1$ as long as
$\gamma$ is in ${\mathcal O}_1$. Consequently
$$
u_\nu(\gamma(t))\ge u_\nu(\gamma(0))-C_1 d_g(\gamma(t), \p
{\mathcal M})= c e^{-u(\gamma(0))} -C_1 d_g(\gamma(t), \p
{\mathcal M})
$$
as along as $\gamma$ is in ${\mathcal O}_1$. Fix an open set
${\mathcal O}_3$ of ${\mathcal M}$ such that $\overline{\mathcal
O}_2\subset {\mathcal O}_3$ and $\overline{\mathcal O}_3\subset
{\mathcal O}_1$. Since $c>0$ and $u\le C_0$ in ${\mathcal O}_1$,
there exist universal constants $0<\delta_1\le\delta_0$ and
$\alpha_0>0$ such that
\begin{equation}\label{4.3}
u_{\nu}\ge \alpha_0 \quad \mbox{in }  {\mathcal O}_3\cap {\mathcal
M}_{\delta_1}.
\end{equation}

Now we are going to introduce a linear elliptic differential
operator ${\mathcal L}_u$ on $\mathcal{M}$. Since $\la_g(U)\in
\Gamma_n$, we can define a tensor $U^{-1}$ on ${\mathcal M}$,
which in local frame has the representation $U^{-1}=\{U^{ij}\}$,
where $U^{ij} U_{jk}=\delta^i_k$ and $\{U_{ij}\}$ denotes the
local representation of $U$. We define
$$
{\mathcal L}_u \psi:=U^{ij} \psi_{ij}-\tr_g(U^{-1})\l \nabla u,
\nabla \psi\r_g, \quad \forall \psi\in C^2({\mathcal M}).
$$
At the end of the proof we will show that for the local function
$u_\nu-c e^{-u}$ there holds
\begin{equation}\label{4.4}
|{\mathcal L}_u (u_\nu-c e^{-u})|\le C\left(1
+\tr_g(U^{-1})\right) \quad  \mbox{in } {\mathcal
M}_{\delta_1}\cap {\mathcal O}_1
\end{equation}
for some universal constant $C$.

We now fix a point $x_0\in {\mathcal O}_2\cap \p {\mathcal M}$ and
consider the function
$$
\psi=A\varphi +B \eta \pm \left(u_\nu-c e^{-u}\right),
$$
where $A$ and $B$ are two large positive constants to be chosen
below, $\varphi(x)=d_g(x, \p {\mathcal M})$, and $\eta:=d_g(x,
x_0)^2$. We can choose a universal constant $0<\delta_2\le
\delta_1$ such that ${\mathcal O}_{\delta_2}(x_0):=\{x\in
{\mathcal M}: d_g(x, x_0)<\delta_2\}\subset {\mathcal O}_3\cap
{\mathcal M}_{\delta_1}$ and $\eta$ is smooth in ${\mathcal
O}_{\delta_2}(x_0)$ with
$$
|\eta|_{C^2({\mathcal O}_{\delta_2}(x_0))}\le C
$$
for some universal constant $C$ independent of $x_0$.

Let us do some calculation first on $\varphi$. Choose a local
orthonormal frame field $\{e_1, \cdots, e_n\}$ around $x_0$ such
that $e_n=\nu$ in ${\mathcal M}_{\delta_2}$. Then $\rho_n=1$ in
${\mathcal M}_{\delta_2}$. Since $\rho$ satisfies the
Hamilton-Jacobi equation $|\nabla \rho|=1$ in ${\mathcal
M}_{\delta_2}$ we have $\rho_i=0$ in ${\mathcal M}_{\delta_2}$ for
$1\le i\le n-1$. Therefore in ${\mathcal M}_{\delta_2}$ there hold
$$
\varphi_n=1 \quad \mbox{and} \quad \varphi_i=0 \quad \mbox{for }
1\le i\le n-1.
$$
This together with (\ref{4.3}) implies that
\begin{equation}\label{4.5}
\l \nabla u, \nabla \varphi\r_g = u_n \ge \alpha_0 \quad \mbox{in
} {\mathcal O}_3\cap {\mathcal M}_{\delta_2}.
\end{equation}
By direct calculation we can see that on $\p {\mathcal M}$ there
holds
$$
\varphi_{ij}=\left\{\begin{array}{lll} -\Pi(e_i, e_j), & 1\le i,
j\le n-1,\\
-h_g, & i=j=n,\\
0, & \mbox{otherwise}.
\end{array}\right.
$$
Since $\p {\mathcal M}$ is totally geodesic in $({\mathcal M},
g)$, we have $\nabla^2 \varphi=0$ on $\p {\mathcal M}$. Therefore,
one may choose a universal constant $0<\delta_3\le \delta_2$ such
that
\begin{equation}\label{4.6}
\nabla^2 \varphi\le \frac{1}{2}\alpha_0 g \quad \mbox{in }
\mathcal{M}_{\delta_3}.
\end{equation}

Returning to the function $\psi$.  Since $u_\nu-c e^{-u}=0$ on $\p
{\mathcal M}$, we have $\psi\ge 0$ on $\p {\mathcal M}\cap
{\mathcal O}_{\delta_3}(x_0)$. Since $\eta\ge \delta_3^2$ on $\p
{\mathcal O}_{\delta_3}(x_0)\cap ({\mathcal M}\backslash \p
{\mathcal M})$ and since $u$ and $|\nabla u|$ are bounded, we also
have $\psi\ge 0 $ on $\p {\mathcal O}_{\delta_3}(x_0)\cap
({\mathcal M}\backslash \p {\mathcal M})$ by choosing $B$ large
enough. Therefore
\begin{equation}\label{4.7}
\psi\ge 0 \quad \mbox{on } \p {\mathcal O}_{\delta_3}(x_0).
\end{equation}

In the following we will show that
\begin{equation}\label{4.8}
{\mathcal L}_u \psi\le 0 \quad \mbox{in } {\mathcal
O}_{\delta_3}(x_0)
\end{equation}
if $A$ is chosen large enough. It is clear that ${\mathcal
L}_u(\eta)\le C\left(1+\tr_g(U^{-1})\right)$. This together with
(\ref{4.4}) gives
\begin{align*}
{\mathcal L}_u\psi&=A{\mathcal L}_u \varphi+B {\mathcal L}_u \eta
\pm {\mathcal L}_u(u_\nu - c e^{-u})\\
& \le A{\mathcal L}_u \varphi +C\left(1+\tr_g(U^{-1})\right).
\end{align*}
But from (\ref{4.5}) and (\ref{4.6}) it follows that
\begin{align*}
{\mathcal L}_u \varphi&= U^{ij}
\varphi_{ij}-\tr_g(U^{-1}) \l\nabla u, \nabla \varphi\r_g \\
&\le \frac{1}{2}\alpha_0 \tr_g (U^{-1})
-\alpha_0 \tr_g(U^{-1})\\
&\le -\frac{1}{2}\alpha_0 \tr_g(U^{-1}).
\end{align*}
Therefore
$$
{\mathcal L}_u\psi\le -\frac{1}{2}A\alpha_0 \tr_g(U^{-1})
+C\left(1+\tr_g(U^{-1})\right) \quad \mbox{in }
\mathcal{O}_{\delta_3}(x_0).
$$
Note that
$$
\tr_g(U^{-1})\ge n\det(g\cdot
U^{-1})^{\frac{1}{n}}=n\det(g^{-1}\cdot U)^{-\frac{1}{n}}=n
e^{2u}\ge \beta_0>0
$$
for some universal constant $\beta_0>0$. Thus
$$
{\mathcal L}_u \psi\le -\frac{1}{2}A\alpha_0 \tr_g(U^{-1})+ C
\tr_g(U^{-1})\le 0 \quad \mbox{in } \mathcal{O}_{\delta_3}(x_0)
$$
if we choose $A$ large enough.

By the maximum principle it follows from (\ref{4.7}) and
(\ref{4.8}) that $\psi\ge 0$ in ${\mathcal O}_{\delta_3}(x_0)$.
Since $\psi(x_0)=0$, we therefore have $\psi_\nu(x_0)\ge 0$.
Consequently $|\nabla^2 u(\nu,\nu)|(x_0)\le C$ for some universal
constant $C$.

In order to complete the proof, we still need to prove claim
(\ref{4.4}). Fix a local orthonormal frame field $\{e_1, \cdots,
e_n\}$ with $e_n=\nu$. For the local function $u_\nu$ it is easy
to check
$$
(u_\nu)_{ij}-u_{nij}=\Ga_{in}^k u_{kj}+\Ga_{jn}^k u_{ki} +b_{ij}^k
u_k,
$$
where
$$
\Ga_{ij}^k =\l \nabla_{e_i} e_j, e_k\r_g \quad \mbox{and}\quad
b_{ij}^k=e_j(\Ga_{in}^k)+\Ga_{jl}^k \Ga_{in}^l
-\Ga_{ji}^l\Ga_{ln}^k,
$$
they are all bounded functions. Recall the commutation formula
$u_{nij}-u_{ijn}=R_{ikjn} u_k$, we therefore have
$$
(u_\nu)_{ij}-u_{ijn}=\Ga_{in}^k u_{kj}+\Ga_{jn}^k u_{ki}
+(b_{ij}^k+R_{ikjn}) u_k
$$
Noting that $(u_\nu)_l=u_{ln}+\Ga_{ln}^k u_k$, we thus have
\begin{align*}
{\mathcal L}_u (u_\nu)&=\left(U^{ij} u_{ijn}-\tr_g(U^{-1}) u_l
u_{ln}\right)+ 2\Ga_{in}^k U^{ij} u_{kj}\\
&\quad +\left(b_{ij}^k+R_{ikjn}\right)U^{ij} u_k
-\tr_g(U^{-1})\Ga_{ln}^k u_l u_k.
\end{align*}
By taking logarithm on equation (\ref{4.1}) and then taking
covariant differentiation, we have
$$
U^{ij}\left(u_{ijn}+2u_{in}u_j-u_l u_{ln} g_{ij}
+(A_g)_{ij,n}\right)=-2n u_n
$$
Consequently
\begin{align*}
{\mathcal L}_u(u_\nu)&=-2n u_n+U^{ij}\left(2\Ga_{in}^k
u_{kj}-2u_{in} u_j-(A_g)_{ij,n}\right)\\
&\quad +\left(b_{ij}^k+R_{ikjn}\right)U^{ij} u_k
-\tr_g(U^{-1})\Ga_{ln}^k u_l u_k.
\end{align*}
Since $(U^{ij})$ is positive definite, we have $|U^{ij}|\le
C\tr_g(U^{-1})$. Using the fact $U^{ij} U_{jk}=\delta^i_k$, the
relation between $U_{ij}$ and $u_{ij}$ and the boundedness of
$|\nabla u|_g$, one can see that
$$
|U^{ij} u_{kj}|+|U^{ij} u_{in}|\le C\left(1+\tr_g(U^{-1})\right)
$$
Therefore
$$
|{\mathcal L}_u(u_\nu)|\le C\left(1+\tr_g(U^{-1})\right) \quad
\mbox{in } {\mathcal M}_{\delta_1}\cap {\mathcal O}_1.
$$
By direct calculation we also have
$$
{\mathcal L}_u(e^{-u})=-e^{-u} U^{ij} u_{ij}+e^{-u}U^{ij} u_i u_j
+e^{-u} \tr_g(U^{-1})|\nabla u|^2
$$
Hence
$$
|{\mathcal L}_u(e^{-u})|\le C\left(1+\tr_g(U^{-1})\right) \mbox{in
} {\mathcal M}_{\delta_1}\cap {\mathcal O}_1.
$$
Putting the above two estimates together, we therefore obtain
(\ref{4.4}).
\end{proof}

Theorem \ref{T4.1} now follows from the combination of Lemma
\ref{L4.1} and the following result which provides more
information than what we really need to complete the proof of
Theorem \ref{T4.1}.

\begin{Lemma}\label{L4.2}
Assume that $(f, \Ga)$ satisfies (\ref{1.1})--(\ref{1.6}) and that
$(\mathcal{M}, g)$ is a smooth compact Riemannian manifold with smooth
umbilic boundary $\p {\mathcal M}$. Let ${\mathcal O}_1$ be an open set of
${\mathcal M}$ and let $u\in C^4({\mathcal O}_1)$ be a solution of
(\ref{1.9.5}). If
\begin{equation}\label{4.9}
|u|\le C_0 \quad \mbox{and} \quad |\nabla u|_g\le C_0 \quad
\mbox{in } \mathcal{O}_1
\end{equation}
and
\begin{equation}\label{4.10}
\nabla^2u(\nu, \nu)\le C_0 \quad \mbox{on } {\mathcal O}_1\cap \p
{\mathcal M}
\end{equation}
for some constant $C_0$, then, for any open set ${\mathcal O}_2$
of ${\mathcal M}$ satisfying $\overline{\mathcal O}_2\subset
{\mathcal O}_1$,
$$
|\nabla^2 u|_g\le C \quad \mbox{in } \mathcal{O}_2
$$
for some constant $C$ depending only on $n$, $c$, $C_0$, $g$, $(f,
\Ga)$, ${\mathcal O}_1$ and ${\mathcal O}_2$.
\end{Lemma}
\begin{proof} Consider the function $w:=e^u$. Recall that
$h_g=0$ on $\p {\mathcal M}$, from (\ref{1.9.5}) it is easy to
check that
\begin{equation}\label{4.11}
\left\{\begin{array}{lll} F(W):=f(\la_g(W))=w^{-1}, \quad
\la_g(W)\in \Ga \quad \mbox{in } {\mathcal O}_1,\\
\\
\frac{\p w}{\p \nu}=c \quad \mbox{on }  {\mathcal O}_1\cap \p
{\mathcal M},
\end{array}\right.
\end{equation}
where
$$
W:=\nabla^2 w-\frac{1}{2w}|\nabla w|_g^2 g+ w A_g.
$$
Moreover, it follows from (\ref{4.9}) and (\ref{4.10}) that
\begin{equation}\label{4.12}
C_1^{-1}\le w\le C_1 \quad \mbox{and} \quad |\nabla w|_g\le C_1
\quad \mbox{in } {\mathcal O}_1
\end{equation}
and
\begin{equation}\label{4.13}
\nabla^2 w(\nu, \nu)\le C_1 \quad \mbox{on } {\mathcal O}_1\cap \p
{\mathcal M}
\end{equation}
for some positive constant $C_1$ depending only on $C_0$. Note
that $\p {\mathcal M}$ is totally geodesic in $({\mathcal M}, g)$.
Thus $\nabla_\tau \nu=0$ on $\p {\mathcal M}$ for any $\tau \in
T(\p {\mathcal M})$. Since $w_\nu=c$ on $\p {\mathcal M}$, for any
$\tau\in T(\p{\mathcal M})$ there holds
$$
\nabla^2 w(\nu, \tau)=\tau(w_\nu)-dw(\nabla_\tau \nu)=0.
$$

Now let $U({\mathcal M})$ denote the unit tangent bundle over
${\mathcal M}$ and consider the function
$$
Q(\xi):=\rho e^{\beta \varphi\circ \pi(\xi)} \nabla^2 w(\xi, \xi),
\quad \xi\in U({\mathcal M}),
$$
where $\pi: T{\mathcal M}\to \mathcal{M}$ is the canonical
projection, $\beta>0$ is a positive constant to be chosen below,
$\varphi\in C^\infty({\mathcal M})$ is a fixed function satisfying
$\varphi(x)=d_g(x, \p {\mathcal M})$ in $\mathcal{M}_{\delta_0}
:=\{x\in {\mathcal M}: d_g(x, \p{\mathcal M})\le \delta_0\}$, and
$\rho\in C_0^\infty({\mathcal O}_1)$ is a cut-off function
satisfying (\ref{2.1.1}) and (\ref{2.1.2}). Suppose the maximum of
$Q$ over $U({\mathcal M})$ is attained at $\bar{\xi}\in
T_{\bar{x}}{\mathcal M}$ for some $\bar{x}\in \mathcal{O}_1$. In
the following we will assume that $\nabla^2 w(\bar{\xi},
\bar{\xi})\ge 1$ since otherwise we are done. We have to consider
two cases: either $\bar{x}\in {\mathcal O}_1\cap \p{\mathcal M}$
or $\bar{x}\in {\mathcal O}_1\backslash \p {\mathcal M}$.

{\it Case 1. $\bar{x}\in {\mathcal O}_1\cap \p {\mathcal M}$}. We
write $\bar{\xi}=\a \nu+\b \tau$, where $\tau$ is a unit vector in
$T_{\bar{x}}(\p {\mathcal M})$ and $\a$ and $\b$ are two numbers
satisfying $\a^2+\b^2=1$. Then by using the maximality of
$\nabla^2 w(\bar{\xi}, \bar{\xi})$ and the fact $\nabla^2w(\nu,
\tau)=0$ one can see that at $\bar{x}$ there holds
\begin{align*}
\nabla^2 w(\bar{\xi}, \bar{\xi})&=\a^2 \nabla^2w(\nu, \nu) +\b^2
\nabla^2 w(\tau, \tau)
+2\a\b \nabla^2 w(\nu, \tau)\\
&\le \left(\a^2+\b^2\right) \nabla^2 w(\bar{\xi}, \bar{\xi})\\
&=\nabla^2w(\bar{\xi}, \bar{\xi}).
\end{align*}
This implies that we can take $\bar{\xi}$ so that either
$\bar{\xi}=\nu$ or $\bar{\xi}\in T(\p {\mathcal M})$.

If $\bar{\xi}=\nu$, then (\ref{4.13}) implies that
$Q(\bar{\xi})\le C$ for some universal constant $C$. So we may
assume that $\bar{\xi}$ is a unit vector in $T_{\bar{x}}(\p
{\mathcal M})$. Choose a local orthonormal frame field $\{e_1,
\cdots, e_n\}$ around $\bar{x}$ so that $e_n=\nu$ on $\p{\mathcal
M}$, and write $\bar{\xi}=\bar{\xi}^i e_i$ at $\bar{x}$. We then
define a vector field $\xi$ near $\bar{x}$ by $\xi=\xi^i e_i$,
where $\xi^i(x)=\bar{\xi}$ for $x$ near $\bar{x}$.  Note that
$\xi^n=0$ near $\bar{x}$. It is clear that $\xi$ is a smooth local
section of $U(\mathcal{M})$. Thus $Q:=Q(\xi)$ has a local maximum
at $\bar{x}$. This implies that
$$
Q_n\le 0 \quad \mbox{at } \bar{x}.
$$
Set $E=\nabla^2 w(\xi, \xi)$. Then, since $\rho_n=0$ and
$\varphi_n=1$ on ${\mathcal O}_1\cap \p {\mathcal M}$,  we have
\begin{equation}\label{4.14}
E_n+\beta E\le 0 \quad \mbox{at } \bar{x}.
\end{equation}
Observe that
\begin{align*}
E_n&=\nabla_n(\xi^i\xi^j w_{ij})=\xi^i\xi^j w_{ijn}+2 \xi^i
\nabla_n \xi^j w_{ij}\\
&=\xi^i\xi^j \left(w_{nij}+R_{kijn} w_k\right) +2\xi^i \nabla_n
\xi^j w_{ij}.
\end{align*}
As calculated in the proof of Lemma \ref{L4.1}, with the same
notations as there we have for $1\le i, j\le n-1$
\begin{align*}
w_{nij}&=(w_\nu)_{ij}+\Ga_{in}^k w_{kj}+\Ga_{jn}^k w_{ki}+
b_{ij}^k
w_k\\
&=\Ga_{in}^k w_{kj}+\Ga_{jn}^k w_{ki}+ b_{ij}^k w_k.
\end{align*}
By the maximality of $Q(\bar{\xi})$ we have $w_{ii}\le  C E$ at
$\bar{x}$ for $1\le i\le n$. Since $\Delta w\ge 0$ in ${\mathcal
M}$, we further have $|\nabla^2 w|\le C E$ at $\bar{x}$. Thus
$E_n\ge -C_3-C_4 E$ for some universal constants $C_3$ and $C_4$.
This together with (\ref{4.14}) implies that
$$
(\beta-C_4) E\le C_3 \quad \mbox{at } \bar{x}.
$$
Therefore $E\le C$ if we choose $\beta>C_4$. Consequently
$Q(\bar{\xi})\le C$.

{\it Case 2. $\bar{x}\in {\mathcal O}_1\backslash \p {\mathcal
M}$.} Choose normal coordinates $x^1, \cdots, x^n$ around
$\bar{x}$ such that
$$
g_{ij}=\delta_{ij} \quad \mbox{and} \quad \frac{\p g_{ij}}{\p
x^k}=0 \quad \mbox{at } \bar{x}.
$$
Moreover, such normal coordinates can be chosen so that
$\{w_{ij}\}$ is diagonal at $\bar{x}$ and $w_{11}=\nabla^2
w(\bar{\xi}, \bar{\xi})$.

Consider the local function $Z:=w_{11}/g_{11}$. By direct
calculation we have
$$
Z_i=w_{11i} \quad \mbox{and} \quad Z_{ij}=w_{11ij} \quad \mbox{at
} \bar{x}.
$$
It is clear that the function
$$
\widetilde{Q}:=\rho e^{\beta \varphi} Z
$$
has a local maximum at $\bar{x}$. Thus at $\bar{x}$ we have
\begin{equation}\label{4.15}
0=\widetilde{Q}_i=\rho e^{\beta\varphi} w_{11i}+\left(\beta
\varphi_i+\frac{\rho_i}{\rho}\right)\widetilde{Q}
\end{equation}
and
\begin{align*}
0&\ge (\widetilde{Q}_{ij})\\
&=\left((\beta\varphi_{ij}-\beta^2
\varphi_i\varphi_j)\widetilde{Q}+\frac{\rho\rho_{ij}-2\rho_i\rho_j}{\rho^2}\widetilde{Q}
-\frac{\rho_i\varphi_j+\rho_j\varphi_i}{\rho} \beta
\widetilde{Q}+\rho e^{\beta \varphi} w_{11ij}\right)
\end{align*}
Let $F^{ij}:=\frac{\p F}{\p W_{ij}}(W)$ and ${\mathcal
T}:=\tr_g(F^{ij})$. Then using $|\nabla \rho|\le C\sqrt{\rho}$ and
the inequality
$$
|w_{11ij}-w_{ij11}|\le C|\nabla^2 w|\le C w_{11}\le \frac{C}{\rho}
\widetilde{Q}.
$$
we have
\begin{equation}\label{4.16}
0\ge e^{-\beta \varphi} F^{ij} \widetilde{Q}_{ij} \ge\rho F^{ij}
w_{11ij}-\frac{C}{\rho} {\mathcal T} \widetilde{Q}\ge \rho F^{ij}
w_{ij11} -\frac{C}{\rho} {\mathcal T} \widetilde{Q}
\end{equation}
By differentiating (\ref{4.11}) twice and using the concavity of
$F$ we get
\begin{align*}
F^{ij}& w_{ij11}-{\mathcal T}\Big(w^{-1} \sum_l w_l w_{l11}+w^{-1}
\sum_l w_{l1}^2-2w^{-2} \sum_l w_l w_{l1} w_1 \\
&+w^{-3} w_1^2 |\nabla w|^2-\frac{1}{2}w^{-2}|\nabla w|^2 w_{11}
\Big) +F^{ij}\left(w_{11} A_{ij} +2w_1 A_{ij,1}+w A_{ij,
11}\right)\\
&\ge -w^{-2}w_{11}+2w^{-3} w_1^2.
\end{align*}
This together with commutation formula and (\ref{4.15}) implies
that
\begin{align*}
\rho F^{ij} w_{ij11}&\ge -C\rho {\mathcal T} w_{11} + \rho
w^{-1}{\mathcal T} \sum_l w_l w_{l11} + \rho w^{-1} {\mathcal T}
w_{11}^2\\
&\ge -C\rho {\mathcal T} w_{11} +\rho w^{-1}{\mathcal T} \sum_l
w_l w_{11l} +\rho w^{-1} {\mathcal T} w_{11}^2\\
&\ge -\frac{C}{\sqrt{\rho}}{\mathcal T} \tilde{Q} +\rho
w^{-1}{\mathcal T} w_{11}^2.
\end{align*}
Thus, it follows from (\ref{4.16}) that $\widetilde{Q}^2\le
C\widetilde{Q}$. Consequently
$Q(\bar{\xi})=\widetilde{Q}(\bar{x})\le C$. The proof is complete.
\end{proof}

\section{\bf Some existence results}
\setcounter{equation}{0}

Let ${\mathcal M}$ be a smooth compact manifold with smooth
boundary $\p {\mathcal M}$. We make use of the double
$\widehat{\mathcal M}$ of ${\mathcal M}$ which is obtained by
gluing two copies of ${\mathcal M}$ along the boundary $\p
{\mathcal M}$. There is a canonical way to make $\widehat{\mathcal
M}$ into a smooth compact manifold without boundary \cite{WRL}.
Given a smooth Riemannian metric $g$ on ${\mathcal M}$, there is a
standard metric $\hat{g}$ on $\widehat{\mathcal M}$ induced from
$g$. In general $\hat{g}$ is only continuous on $\widehat{\mathcal
M}$. However, if $\p {\mathcal M}$ is totally geodesic in
$({\mathcal M}, g)$, then $\hat{g}$ is $C^{2,1}$ on
$\widehat{\mathcal M}$, see \cite[Appendix]{E} for instance.

\subsection{Proof of Theorem \ref{T1.1}}

We may assume that $({\mathcal M}, g)$ is not conformally
equivalent to the standard half sphere ${\mathbb S}^n_+$ since
otherwise the existence result is obvious.

First note that we may assume $\la(A_g)\in \Ga$ on ${\mathcal M}$
and $h_g>0$ on $\p {\mathcal M}$ in the following argument. To see
this, consider the metric $g_\varepsilon :=(1-\varepsilon
\varphi)^{\frac{4}{n-2}}g$, where $\varepsilon>0$ is a small
number, and $\varphi\in C^\infty({\mathcal M})$ is a function such
that $\varphi(x)=d_g(x, \p {\mathcal M})$ when $d_g(x, \p
{\mathcal M})\le \delta_0$. Since $\la(A_g)\in \Ga$ on ${\mathcal
M}$, we can fix an $\varepsilon>0$ small enough so that
$\la(A_{g_\varepsilon})\in \Ga$ on ${\mathcal M}$. Then noting
that $\varphi=0$, $\frac{\p \varphi}{\p \nu}=1$ and $h_g\ge 0$ on
$\p {\mathcal M}$ we have
$$
h_{g_\varepsilon}=-\frac{2}{n-2}\frac{\p}{\p
\nu}(1-\varepsilon\varphi)+h_g=\frac{2\varepsilon}{n-2}+h_g>0
\quad \mbox{on } \p {\mathcal M},
$$

Since $R_g>0$ on ${\mathcal M}$ and $h_g>0$ on $\p {\mathcal M}$,
one can find a metric $g_0$ conformal to $g$ such that $R_{g_0}>0$
on ${\mathcal M}$ and $h_{g_0}=0$ on $\p {\mathcal M}$, see
\cite[Theorem 0.1]{HL} for instance. Write $g=e^{-2\varphi} g_0$
for some function $\varphi\in C^\infty({\mathcal M})$. Let $u\in
C^\infty({\mathcal M})$ be a solution of (\ref{1.9}) with $c=0$
and let $\tilde{g}:=e^{-2u} g$. Then $\tilde{g}=e^{-2v}g_0$ with
$v=u+\varphi$ and $h_{\tilde{g}}=0$ on $\p {\mathcal M}$. Let
$\widehat{\mathcal M}$ denote the double of ${\mathcal M}$, and
let $\hat{g}_0$ and $\hat{\tilde{g}}$ denote the standard metrics
on $\widehat{\mathcal M}$ induced from $g_0$ and $\tilde{g}$
respectively. Since $\p {\mathcal M}$ is totally geodesic in both
$({\mathcal M}, g_0)$ and $({\mathcal M}, \tilde{g})$, $\hat{g}_0$
and $\hat{\tilde{g}}$ are in $C^{2,1}(\widehat{\mathcal M})$.
Moreover, $\hat{\tilde{g}}$ is still conformal to $\hat{g}_0$ with
$\hat{\tilde{g}}=e^{-2\hat{v}} \hat{g}_0$ for some function
$\hat{v}\in C^{2, 1}(\widehat{\mathcal M})$, and
$$
f(\la(A_{e^{-2\hat{v}}\hat{g}_0}))=1, \quad
\la(A_{e^{-2\hat{v}}\hat{g}_0})\in \Ga \quad \mbox{on }
\widehat{\mathcal M}.
$$
Since $({\mathcal M}, g_0)$ is locally conformally flat, so is
$(\widehat{\mathcal M}, \hat{g}_0)$. Note that $R_{\hat{g}_0}>0$
on $\widehat{\mathcal M}$ and $(\widehat{\mathcal M}, \hat{g}_0)$
is not conformally equivalent to ${\mathbb S}^n$. Therefore, it
follows from the proof of \cite[Theorem 1]{LL03b} that
$$
|\nabla \hat{v}|\le C\quad \mbox{and} \quad \hat{v}\ge -C\quad
\mbox{on } \widehat{\mathcal M}
$$
for some universal constant $C$. However, an upper bound for
$\hat{v}$ is not yet available since we do not have
$\la(A_{\hat{g}_0})\in \Ga$. Since $\hat{v}=v$ and $v=u+\varphi$
on ${\mathcal M}$, we have
\begin{equation}\label{5.1}
|\nabla u|\le C_0\quad \mbox{and} \quad u\ge -C_0\quad  \mbox{on }
{\mathcal M}.
\end{equation}
for some universal constant $C_0$.

Returning to problem (\ref{1.9}) with $c=0$. Suppose the minimum
of $u$ over ${\mathcal M}$ is attained at some point $x_0\in
{\mathcal M}$. Since $h_g>0$ on $\p {\mathcal M}$, we have $x_0\in
{\mathcal M}\backslash \p {\mathcal M}$. Thus $\nabla u=0$ and
$\nabla^2 u\ge 0$ at $x_0$. Since $\la(A_g)\in \Ga$, we therefore
have $e^{-2u} \ge f(\la(A_g))>0$. Consequently there is a
universal constant $C$ such that
\begin{equation}\label{5.1.5}
\min_{\mathcal M} u\le C
\end{equation}
Combining (\ref{5.1}) and (\ref{5.1.5}) gives
$$
-C\le u\le C \quad \mbox{and} \quad |\nabla u|\le C \quad \mbox{on
} {\mathcal M}.
$$
Therefore it follows from Remark \ref{R3.1} that
$$
\|u\|_{C^{4,\a}({\mathcal M})}\le C
$$
for some universal constant $C$.

Now we will use the degree theory argument to prove the existence.
To this end, as in \cite{LL03}, for each $0\le t\le 1$ let
$$
f_t(\la):=f(t\la+(1-t)\sigma_1(\la) e)
$$
which is defined on
$$
\Ga_t:=\left\{\la\in {\mathbb R}^n: t\la+(1-t)\sigma_1(\la) e\in
\Ga\right\},
$$
where $e=(1,1, \cdots, 1)$.

We now consider the problem
$$
\left\{\begin{array}{lll} f_t(\la(A_{g_u}))=1, \quad
\la(A_{g_u})\in \Ga_t \quad \mbox{on } {\mathcal M},\\
\\
h_{g_u}=0 \quad \mbox{on } {\mathcal M},
\end{array}\right.
$$
where $g_u=e^{-2u} g$ for some smooth function $u$ on ${\mathcal
M}$. From the above argument we have already obtained
$\|u\|_{C^{4, \a}({\mathcal M})}\le C$ for some universal constant
independent of $t$. Now we  set
\begin{align*}
{\mathcal O}_t^*=\Big\{&u\in C^{4,\a}({\mathcal
M}):\la(A_{g_u})\in \Ga_t, \|u\|_{C^{4,\a}({\mathcal M})}<2C,\\
& \frac{1}{2}<f_t(\la(A_{g_u}))<2 \mbox{ on } {\mathcal M} \mbox{
and } \frac{\p u}{\p \nu}=0 \mbox{ on } \p {\mathcal M}\Big\}.
\end{align*}
Define $F_t: {\mathcal O}^*_t\to C^{2,\a}({\mathcal M})$ by
$F_t[u]:=f_t(\la(A_{g_u}))-1$. It follows from \cite{Li89} that
$\deg(F_t, {\mathcal O}^*_t, 0)$ is well-defined and is
independent of $t$. But when $t=0$ the corresponding problem is
the Yamabe problem with boundary. Based on \cite{Sch} it was shown
in \cite{HL} that $\deg(F_0, {\mathcal O}^*_0, 0)=-1$. Therefore
$\deg(F_1, {\mathcal O}^*_1, 0)=-1\ne 0$. The proof is thus
complete.

\subsection{Proof of Theorem \ref{T1.2}}

The proof of Theorem \ref{T1.2} is based on some lemmas in the following.

\begin{Lemma}\label{L5.1}
Let $(f, \Ga)$ satisfy (\ref{1.1})--(\ref{1.6}) and let
$({\mathcal M}, g)$ be a smooth compact Riemannian manifold with smooth
boundary $\p {\mathcal M}$. Suppose $\p{\mathcal M}$ is umbilic and $
({\mathcal M}, g)$ is locally conformally flat near $\p {\mathcal M}$.
Let $u\in C^4({\mathcal M})$ be a solution of (\ref{1.11}) with
$\eta$ being positive. If $|u|\le C_0$ on ${\mathcal M}$, then
$$
|\nabla u|+|\nabla^2 u|\le C \quad \mbox{on } {\mathcal M}
$$
for some constant $C$ depending only on $n$, $C_0$, $g$, $\psi$,
$\eta$, and $(f, \Ga)$.
\end{Lemma}
\begin{proof}
This is the combination of Theorem \ref{T1.3} and Theorem
\ref{T1.5}.
\end{proof}

\begin{Lemma}\label{L5.2}
Let $(f, \Ga)$ and $({\mathcal M}, g)$ be as in Lemma \ref{L5.1}
with $\la(A_g)\in \Ga$ on ${\mathcal M}$ and $h_g\ge 0$ on $\p
{\mathcal M}$. Then problem (\ref{1.11}) with $\psi(x,
z)=\psi_0(x) e^{az}$ and $\eta(x, z)=\eta_0(x) e^{bz}$ has a
unique solution, where $a$ and $b$ are positive constants, and
$\psi_0\in C^2({\mathcal M})$ and $\eta_0\in C^2(\p {\mathcal M})$
are positive functions.
\end{Lemma}
\begin{proof} By perturbing $g$ as in the proof of Theorem
\ref{T1.1}, we may assume that $\la(A_g)\in \Ga$ on ${\mathcal M}$
and $h_g>0$ on $\p {\mathcal M}$. From the maximum principle, it
is easy to check that there is a positive universal constant $C$
such that $-C\le u\le C$ on ${\mathcal M}$ for any solution $u$ of
(\ref{1.11}) with $\psi(x, z)=\psi_0(x) e^{az}$ and $\eta(x,
z)=\eta_0(x) e^{bz}$. Therefore, it follows from Lemma \ref{L5.1}
and  the result of Lieberman-Trudinger \cite{LT86} that we have
uniform $C^{2,\a}({\mathcal M})$ estimates on $u$. Since $a>0$,
$b>0$, $\psi_0$ and $\eta_0$ are positive, the linearized problem
is uniquely solvable. Therefore, the method of continuity
concludes the existence and uniqueness.
\end{proof}

Next we will use the recent results of Trudinger-Wang in
\cite{TW05} to establish a Harnack type inequality.

\begin{Lemma}\label{L5.3}
Let $({\mathcal M}, g)$ be a smooth compact Riemannian manifold with
smooth boundary $\p {\mathcal M}$. Suppose $\p {\mathcal M}$ is umbilic
and $({\mathcal M}, g)$ is locally conformally flat near 
$\p {\mathcal M}$. For $k>\frac{n}{2}$, let $[g]_k^+$ 
denote the set of $C^\infty$
metrics $\tilde{g}$ conformal to $g$ such that
$\la(A_{\tilde{g}})\in \Ga_k$ on ${\mathcal M}$ and
$h_{\tilde{g}}\ge 0$ on $\p {\mathcal M}$. If $({\mathcal M}, g)$
is not conformally equivalent to the standard half sphere
${\mathbb S}^n_+$, then there is a positive constant $C$ depending
only on $k$ and $({\mathcal M}, g)$ such that for any metric
$\tilde{g}:=\chi g\in [g]_k^+$ there holds
\begin{equation}\label{5.3}
\max_{\mathcal M} \chi\le C\min_{\mathcal M} \chi.
\end{equation}
\end{Lemma}

\begin{proof}
Let $[g]_k^*$ denote the set of metrics $\tilde{g}$ conformal to
$g$ such that $\la(A_{\tilde{g}})\in \Ga$ on ${\mathcal M}$ and
$h_{\tilde{g}}>0$ on $\p {\mathcal M}$. We remark that it suffices
to establish the Harnack inequality (\ref{5.3}) for
$\tilde{g}=\chi g\in [g]_k^*$. Indeed, for any metric
$\tilde{g}=\chi g\in [g]_k^+$, as in the proof of Theorem
\ref{T1.1} we can find a function $\varphi\in C^\infty({\mathcal
M})$ with $\frac{1}{2}\le \varphi\le 1$ on ${\mathcal M}$ such
that $\varphi \tilde{g}=(\varphi \chi) g\in [g]_k^*$. Thus we have
$$
\frac{1}{2}\max_{\mathcal M}\chi\le \max_{\mathcal M} (\varphi
\chi)\le C \max_{\mathcal M} (\varphi \chi) \le C\min_{\mathcal M}
\chi,
$$
which gives the desired inequality.

By a conformal deformation of $g$ without loss of generality we
may assume that $h_g=0$ on $\p {\mathcal M}$. Since $\p {\mathcal
M}$ is umbilic, it must be totally geodesic in $({\mathcal M},
g)$. Let $\widehat{\mathcal M}$ be the double of ${\mathcal M}$.
For any metric $\tilde{g}$ on ${\mathcal M}$, there is a standard
metric $\hat{\tilde{g}}$ on $\widehat{\mathcal M}$ induced by
$\tilde{g}$. In general $\hat{\tilde{g}}$ is only continuous on
$\widehat{\mathcal M}$. However, since $\p {\mathcal M}$ is
totally geodesic in $({\mathcal M}, g)$, it follows from
\cite[Appendix]{E} that $\hat{g}$ is $C^{2,1}$ on
$\widehat{\mathcal M}$.

For any metric $\tilde{g}\in [g]_k^*$, note that $\hat{\tilde{g}}$
is conformal to $\hat{g}$, we may write
$\hat{\tilde{g}}=e^{-2\hat{w}} \hat{g}$ for some function
$\hat{w}\in C^0(\widehat{\mathcal M})$ which is $C^\infty$ in
$\widehat{\mathcal M}\backslash \p {\mathcal M}$. Let
$w_1:=\hat{w}|_{\mathcal M}$ and $w_2:=w|_{\widehat{\mathcal
M}\backslash({\mathcal M}\backslash \p {\mathcal M})}$. Since
$h_{\tilde{g}}> 0$ on $\p {\mathcal M}$, we have $\frac{\p w_1}{\p
\nu_1}> 0$ and $\frac{\p w_2}{\p \nu_2}> 0$ on $\p {\mathcal M}$,
where $\nu_1$ and $\nu_2$ denote the inward unit normal vector
fields to $\p {\mathcal M}$ in $({\mathcal M}, g)$ and
$(\widehat{\mathcal M}\backslash ({\mathcal M}\backslash \p
{\mathcal M}), \hat{g})$ respectively. Thus we may smoothly extend
$w_1$ and $w_2$ to a neighborhood ${\mathcal U}$ of $\p {\mathcal
M}$ in $\widehat{\mathcal M}$ so that $w_i\le w$ and
$\la(A_{\hat{g}_{w_i}})\in \Ga_k$ on ${\mathcal U}$ for $i=1, 2$,
where $\hat{g}_{w_i}:=e^{-2w_i} \hat{g}$. Therefore, noting that
our background metric $\hat{g}$ is $C^{2,1}$ on $\widehat{\mathcal
M}$ and $\hat{w}=\max\{w_1, w_2\}$ on ${\mathcal U}$, we can apply
\cite[Lemma 3.7]{TW05} to conclude that
$\hat{\tilde{g}}=e^{-2\hat{w}}\hat{g}$ is $k$-admissible in the
sense of Trudinger-Wang (This will be simply called $k$-admissible
in the sequel).

In the following we will follow the idea in \cite{TW05} to give
the proof of Lemma \ref{L5.3}. Suppose the Harnack inequality
(\ref{5.3}) does not hold. Then there is a sequence of smooth
metrics $g_j:=e^{-2w_j} g\in [g]_k^*$ such that
$$
\max_{\mathcal M} w_j-\min_{\mathcal M} w_j\ge j, \quad j=1, 2,
\cdots.
$$
By subtracting a constant if necessary, we may assume
$\max_{\mathcal M} w_j=0$. Then $\min_{\mathcal M} w_j\rightarrow
-\infty$ as $j\rightarrow \infty$. Consider the function
$\hat{w}_j$ on $\widehat{\mathcal M}$ induced by $w_j$ through
$\hat{g}_j=e^{-2\hat{w}_j} \hat{g}$. From the above argument we
know $\hat{g}_j$ is $k$-admissible. Thus \cite[Lemma 3.1]{TW05}
shows that for $\hat{u}_j:=e^{\hat{w}_j}$ there holds the
H\"{o}lder estimate
$$
\frac{|\hat{u}_j(x)-\hat{u}_j(y)|}{d_{\hat{g}}(x,y)^\a}\le
C\int_{\widehat{\mathcal M}} \hat{u}_j d\mu_{\hat{g}}
$$
for some $0<\a\le 2-\frac{n}{k}$,  where $\a$ and $C$ are
independent of $j$. Note that $0<\hat{u}_j\le 1$. By Arzela-Ascoli
theorem we may assume that $\hat{u}_j\rightarrow \hat{u}$
uniformly on $\widehat{\mathcal M}$ for some function $\hat{u}\in
C^\a(\widehat{\mathcal M})$ with $0\le \hat{u}\le 1$. Since
$\max_{\mathcal M}\hat{u}_j=1$, we have $\hat{u}\ne 0$. Define
$\hat{w}:=\log \hat{u}$, and let
$$
S_{\hat{w}}:=\bigcap_{\beta<0}\left\{x\in \widehat{\mathcal M}:
\hat{w}(x)<-\beta\right\}
$$
which is called the set of singularity points of $\hat{w}$. Since
$\min_{\mathcal M} \hat{u}_j\rightarrow -\infty$ as $j\rightarrow
\infty$, we know $S_{\hat{w}}\ne \emptyset$. Moreover, since each
$\hat{g}_j$ is $k$-admissible, as the limit
$\hat{g}_{\hat{w}}:=e^{-2\hat{w}} \hat{g}$ is also $k$-admissible
on $\widehat{\mathcal M}$.

Since $\hat{g}$ is smooth away from $\p {\mathcal M}$ and is
locally conformally flat near $\p {\mathcal M}$, by the
$k$-admissibility of $\hat{g}_{\hat{w}}$, the argument in
\cite{TW05} shows that near any singularity point $x_0$ of
$\hat{w}$ there holds
\begin{equation}\label{5.4}
\hat{w}(x)=2\log |x-x_0|+o(1)
\end{equation}
in a normal neighborhood of $x_0$; moreover, the singularity points
are isolated. For a fixed point $y\in \widehat{\mathcal
M}\backslash S_{\hat{w}}$, by using the Bishop volume comparison
theorem and an approximation argument it was shown in \cite[Lemma
3.4]{TW05} that the ratio
$$
Q(r):=\frac{\mbox{Vol}_{\hat{g}_{\hat{w}}}\left(
B_{y,r}[\hat{g}_{\hat{w}}]\right)}{r^n}\le \omega_n, \quad
0<r<\infty,
$$
where $B_{y,r}[\hat{g}_{\hat{w}}]$ denotes the geodesic ball in
${\mathcal M}$ of radius $r$ with center at $y$, and $\omega_n$ is
the volume of the unit ball in ${\mathbb R}^n$. But by (\ref{5.4}) it
was shown in \cite[Lemma 3.4]{TW05} that each singularity point of
$\hat{w}$ contributes a factor $\omega_n$ to the ratio $Q(r)$.
Therefore $S_{\hat{w}}$ must consists of a single point, say
$S_{\hat{w}}=\{x_0\}$, and $Q(r)\equiv\omega_n$. Moreover, noting
that the symmetry of $\hat{w}$ with respect to $\p {\mathcal M}$,
we must have $x_0\in \p {\mathcal M}$.

Next we are going to show that $\hat{w}$ is $C^{2,1}$ away from
$x_0$. Since $\hat{g}$ is smooth in $\widehat{\mathcal
M}\backslash \p {\mathcal M}$, the argument of \cite[Lemma
3.5]{TW05} can be applied directly to show that $\hat{w}\in
C^{1,1}(\widehat{\mathcal M}\backslash \p {\mathcal M})$. However,
since $\hat{g}$ is only $C^{2,1}$ across $\p {\mathcal M}$, when
we consider the regularity at a point $y_0\in \p {\mathcal
M}\backslash \{x_0\}$, we need to check carefully the proof of
\cite[Lemma 3.5]{TW05} when using the existence result on a
Dirichlet problem in \cite{G05}. We may choose a neighborhood
${\mathcal O}_1$ of $y_0$ in $\widehat{\mathcal M}\backslash
\{x_0\}$ on which $\hat{g}$ is conformally flat, i.e.
$\hat{g}=e^{-2\eta} \hat{g}_0$ for some function $\eta\in
C^{2,1}({\mathcal O}_1)$, where $\hat{g}_0$ is the flat metric.
Then $\hat{g}_{\hat{w}}=e^{-2\hat{v}} \hat{g}_0$ with
$\hat{v}=\hat{w}+\eta$. Since $\hat{g}_{\hat{w}}$ is
$k$-admissible, there is a sequence of $k$-admissible metrics
$\hat{g}_j:=e^{-2{\hat{v}_j}} \hat{g}_0$ with $\hat{v_j}$ smooth
on ${\mathcal O}_1$ and $\hat{v}_j\rightarrow \hat{v}$ uniformly
on ${\mathcal O}_2$ for some neighborhood ${\mathcal O}_2$ of
$y_0$ satisfying $\overline{\mathcal O}_2\subset {\mathcal O}_1$.
Let $\{\varepsilon_j\}$ be a sequence of positive numbers such
that $\varepsilon_j\searrow 0$ and
$$
0<\varepsilon_j <\sigma_k\left(\la_{\hat{g}_0} \left(\nabla_0^2
\hat{v}_j+d\hat{v}_j\otimes d \hat{v}_j
-\frac{1}{2}|\nabla_0\hat{v}_j|^2 \hat{g}_0\right)\right).
$$
where $\nabla_0$ denotes the Levi-Civita connection of
$\hat{g}_0$. Consider the problem
$$
\left\{\begin{array}{lll}
\sigma_k\left(\la_{\hat{g}_0}\left(\nabla_0^2\hat{\varphi}_j
+d\hat{\varphi}_j\otimes
d \hat{\varphi}_j-\frac{1}{2} |\nabla_0\hat{\varphi}_j|^2
\hat{g}_0\right)\right)=\varepsilon_j \quad \mbox{on } {\mathcal
O}_2\\
\\
\hat{\varphi}_j=\hat{v}_j \quad \mbox{on } \p {\mathcal O}_2.
\end{array}\right.
$$
It follows from \cite{G05} that such $\hat{\varphi}_j$ exists,
$\|\hat{\varphi}_j\|_{C^2({\mathcal O}_2)}\le C$ and
$\{\hat{\varphi}_j\}$ is monotone increasing. Define
$\hat{\varphi}:=\lim_{j\rightarrow \infty} \hat{\varphi}_j$, then
$\hat{\varphi}\in C^{1,1}({\mathcal O}_2)$. As shown in
\cite[Lemma 3.5]{TW05} $\hat{\varphi}=\hat{v}$ on ${\mathcal
O}_2$. Therefore $\hat{w}=\hat{v}-\eta\in C^{1,1}({\mathcal
O}_2)$. Combining the above we obtain $\hat{w}\in
C^{1,1}(\widehat{\mathcal M}\backslash\{x_0\})$.  By the symmetry
of $\hat{w}$ with respect to $\p {\mathcal M}$ we have $\frac{\p
\hat{w}}{\p \nu}=0$ on $\p {\mathcal M}\backslash \{x_0\}$.
Moreover the argument in \cite[Lemma 3.5]{TW05} gives
$R_{\hat{g}_{\hat{w}}}=0$ a.e. on $\widehat{\mathcal M}$, where
$R_{\hat{g}_{\hat{w}}}$ denotes the scalar curvature of
$\hat{g}_{\hat{w}}$. Thus by the Yamabe equation we know
$\hat{v}:=e^{-\frac{n-2}{2}\hat{w}}$ satisfies
$$
\left\{\begin{array}{lll} -\Delta_g \hat{v}+\frac{n-2}{4(n-1)}
R_{\hat{g}} \hat{v}=0 \quad \mbox{in
} {\mathcal M}\backslash \{x_0\},\\
\\
\frac{\p \hat{v}}{\p \nu}=0 \quad \mbox{on } \p {\mathcal
M}\backslash \{x_0\}.
\end{array}\right.
$$
Note that $\hat{g}$ and $R_{\hat{g}}$ are smooth on ${\mathcal
M}$, it follows from the regularity theory theory for uniformly
elliptic equation with Neumann boundary condition that $\hat{v}\in
C^\infty({\mathcal M}\backslash\{x_0\})$. Thus $\hat{w}\in
C^\infty({\mathcal M}\backslash \{x_0\})$. By symmetry we
correspondingly have $\hat{w}\in C^\infty((\widehat{\mathcal
M}\backslash ({\mathcal M}\backslash\p {\mathcal M}))\backslash
\{x_0\})$. Therefore $\hat{w}\in C^{2,1}(\widehat{\mathcal
M}\backslash\{x_0\})$.

By the $k$-admissibility of $\hat{g}_{\hat{w}}$ and the result in
\cite{GVW03} we know $\hat{g}_{\hat{w}}$ has nonnegative Ricci
curvature. The asymptotic formula (\ref{5.4}) implies that
$(\widehat{\mathcal M}\backslash \{x_0\}, \hat{g}_{\hat{w}})$ is a
complete manifold with $Q(r)=\omega_n$. Hence $(\widehat{\mathcal
M}\backslash\{x_0\}, \hat{g}_{\hat{w}})$ is isometric to the
Euclidean space. Consequently $({\widehat{\mathcal M}}, \hat{g})$
is conformally equivalent to ${\mathbb S}^n$. This contradicts the
assumption that $({\mathcal M}, g)$ is not conformally equivalent
to ${\mathbb S}^n_+$.
\end{proof}

\begin{proof}[Proof of Theorem \ref{T1.2}]
Without loss of generality, we may assume that $({\mathcal M}, g)$
is not conformally equivalent to ${\mathbb S}^n_+$. As indicated
in the proof of Theorem \ref{T1.2}, we may assume $\la(A_g)\in
\Ga$ on ${\mathcal M}$ and $h_g>0$ on $\p {\mathcal M}$. We will
adapt the idea in the proof of \cite[Theorem C]{TW05} to complete
the argument.

For a positive function $v\in C^2({\mathcal M})$ we will use the
notation
\begin{equation}\label{5.10}
 V[v]=-\nabla^2 v+\frac{n}{n-2} \frac{\nabla v\otimes
\nabla v}{v}-\frac{1}{n-2}\frac{|\nabla v|^2}{v} g +\frac{n-2}{2}
v A_g.
\end{equation}
Note that for the metric $g_v:=v^{\frac{4}{n-2}} g$ we have
$$
A_{g_v}=\frac{2}{n-2}v^{-1} V[v] \mbox{ on } {\mathcal M}  \mbox{
and } h_{g_v}=v^{-\frac{n}{n-2}}\left(-\frac{2}{n-2}\frac{\p
v}{\p\nu}+h_g v\right) \mbox{ on } \p {\mathcal M}.
$$
Thus, if we can prove the existence of a positive function $v\in
C^2({\mathcal M})$ such that
\begin{equation}\label{5.5}
\left\{\begin{array}{lll} f(\la(V[v]))=\frac{n-2}{2}\varphi_0 v^p,
\quad \la(V[v])\in \Ga
\quad \mbox{on } {\mathcal M},\\
\\
\frac{\p v}{\p \nu}-\frac{n-2}{2}h_g v=-\frac{n-2}{2}h_0 v^q
\quad \mbox{on } \p {\mathcal M},
\end{array}\right.
\end{equation}
where $p=\frac{n+2}{n-2}$, $q=\frac{n}{n-2}$, and $\la(V[v])$
denote the eigenvalues of $V[v]$ with respect to $g$, then
$u=-\frac{2}{n-2}\log v$ is a solution of (\ref{5.2}).

Since $\la(A_g)\in \Ga$ on ${\mathcal M}$, $h_g>0$ on $\p
{\mathcal M}$, $p>1$ and $q>1$, it is always possible to find
positive numbers $\a_0$ and $\varepsilon_0$ such that
\begin{equation}\label{5.6}
\left\{\begin{array}{lll} f(\la(V[\a_0]))>\frac{n-2}{2}\varphi_0
\left(\a_0^2+\varepsilon_0^2\right)^{p/2} \quad \mbox{on }
{\mathcal
M},\\
\\
- h_g \a_0 <-h_0(\a_0^2 +\varepsilon_0^2)^{q/2} \quad \mbox{on }
\p {\mathcal M}.
\end{array}\right.
\end{equation}
Let $\varepsilon\in (0,\varepsilon_0]$ be any number. We
consider the auxiliary problem ($P_{t, \varepsilon}$) as follows
$$
\left\{\begin{array}{lll} f(\la(V[v]))=\frac{n-2}{2} t
\varphi_0\left(v^2+\varepsilon^2\right)^{p/2}, \quad  \la(V[v])\in
\Ga \quad
\mbox{on } {\mathcal M},\\
\\
\frac{\p v}{\p \nu}-\frac{n-2}{2}h_g v=-\frac{n-2}{2} t h_0
\left(v^2+\varepsilon^2 \right)^{q/2} \quad \mbox{on } \p
{\mathcal M},
\end{array}\right.
$$
where $t>0$ is a parameter.

\vskip 0.2cm

{\it Claim 1}. For any $t_0>0$ there exists a positive constant
$C$ independent of $\varepsilon$ such that any solution
$v$ of ($P_{t,\varepsilon}$) with $t\ge t_0$
and $0<\varepsilon\le \varepsilon_0$ satisfies $v
\le C$ on ${\mathcal M}$. Moreover, for each
$0<\varepsilon\le \varepsilon_0$ there exists $\bar{t}>1$ such that
($P_{t, \varepsilon}$) has no solution for $t\ge \bar{t}$.

\vskip 0.2cm

Indeed, suppose there exist two sequences $\{t_j\}$ and
$\{\varepsilon_j\}$ satisfying $t_j\ge t_0$ and $0<\varepsilon_j\le
\varepsilon_0$ and a solution $v_j$ of ($P_{t_j, \varepsilon_j}$)
such that $\sup_{\mathcal M} v_j\rightarrow \infty$. Then by Lemma \ref{L5.3} we have
$m_j:=\inf_{\mathcal M} v_j\rightarrow \infty$. Note that for the
function $\tilde{v}_j:=v_j/m_j$
$$
\left\{\begin{array}{lll} f(\la(V[\tilde{v_j}]))\ge \frac{n-2}{2}
t_0 \varphi_0  m_j^{p-1}\rightarrow
\infty \quad \mbox{ on } {\mathcal M},\\
\\
\frac{\p \tilde{v}_j}{\p \nu}-\frac{n-2}{2}h_g \tilde{v}_j \le
-\frac{n-2}{2} t_0 h_0 m_j^{q-1}\rightarrow -\infty  \quad \mbox{ on }
\p {\mathcal M}.
\end{array}\right.
$$
By Lemma \ref{L6.1} we have $\inf_{\mathcal M}
\tilde{v}_j\rightarrow \infty$. This is a contradiction since
$\inf_{\mathcal M} \tilde{v}_j=1$. For the second assertion, let
$t>1$ and let $v$ be a solution of ($P_{t, \varepsilon}$). Then
$$
\left\{\begin{array}{lll} f(\la(V[v]))\ge \frac{n-2}{2}
t\varphi_0 \varepsilon^p \quad \mbox{on } {\mathcal M},\\
\\
\frac{\p v}{\p \nu}-\frac{n-2}{2} h_g v
\le -\frac{n-2}{2} t h_0 \varepsilon^q \quad \mbox{on } \p {\mathcal M}.
\end{array}\right.
$$
By Lemma \ref{L6.1} again this implies that $v\ge c_0 t$
for some positive constant $c_0$ independent of $t$. Thus
($P_{t, \varepsilon}$) has no solution if $t$ is large enough
since $v_t$ is uniformly bounded from above.

\vskip 0.2cm

It is important to note that the constant $C$ in Claim 1 does not depend on
$\varepsilon$. Unless stated otherwise, constants appeared below 
allow $\varepsilon$-dependence. From now on we denote ($P_{t,
\varepsilon}$) simply by ($P_t$).

We now define the mapping $T_t: C^2({\mathcal M})\to C^2({\mathcal
M})$ so that for any $v_1\in C^2({\mathcal M})$, $T_t(v_1)$ is the
solution of
$$
\left\{\begin{array}{lll} f(\la(V[v]))=\frac{n-2}{2} t \varphi_0
\left(v_1^2+\varepsilon^2\right)^{p/2} \quad \la(V[v])\in \Ga
\quad \mbox{on }
{\mathcal M},\\
\\
\frac{\p v}{\p \nu}-\frac{n-2}{2}h_g v=-\frac{n-2}{2} t h_0
\left(v_1^2+\varepsilon^2 \right)^{q/2} \quad \mbox{on } \p
{\mathcal M}.
\end{array}\right.
$$
From Lemma \ref{L5.2} it follows that $T_t$ is well-defined.
Moreover $T_t(v_1)$ is a positive function. By a priori estimates
one can see that $T_t$ is a compact operator for each $t>0$. Note
that $T_t=t T_1$ for $t>0$, we may continuously extend $T_t$ to
$t=0$ by setting $T_0=0$.

For the number $\a_0>0$ satisfying (\ref{5.6}) we set
$$
\Phi:=\left\{v\in C^2({\mathcal M}): |v|<\a_0 \mbox{ on }
{\mathcal M}\right\}.
$$

\vskip 0.2cm

{\it Claim 2}. There exists a large number $R_0$ independent of $t$
such that for any $R\ge R_0$
$$
(I-T_t)^{-1}(0)\cap \p (\Phi\cap B_R)=\emptyset
$$
for $0\le t\le 1$, where $B_R:=\{\varphi\in C^2({\mathcal M}):
\|\varphi\|_{C^2({\mathcal M})}<R\}$.

\vskip 0.2cm

To see this, let $v\in (I-T_t)^{-1}(0)\cap \p (\Phi\cap B_R)$ for some
$0<t\le 1$. This implies that $v$ is a
solution of ($P_t$) and $0<v\le \a_0$ on ${\mathcal M}$. By using
(\ref{5.6}) we have
$$
\left\{\begin{array}{lll} f(\la(V[v]))\le \frac{n-2}{2} \varphi_0
\left(\a_0^2+\varepsilon^2\right)^{p/2}< f(\la(V[\a_0])) \quad
\mbox{on } {\mathcal M},\\
\\
\frac{\p v}{\p \nu}-\frac{n-2}{2}h_g v\ge -\frac{n-2}{2} h_0
\left(\a_0^2+\varepsilon^2\right)^{q/2}>\frac{\p \a_0}{\p
\nu}-\frac{n-2}{2}h_g \a_0 \quad \mbox{on } \p {\mathcal M}.
\end{array}\right.
$$
Therefore Lemma \ref{L6.1} implies that $0<v<\a_0$. Consequently
$v\in \p B_R$, i.e.
\begin{equation}\label{5.7}
\|v\|_{C^2({\mathcal M})}=R.
\end{equation}
However, note that the function $\tilde{v}:=t^{-1} v$ satisfies
$$
\left\{\begin{array}{lll} \frac{n-2}{2}\varphi_0 \varepsilon^p \le
f(\la(V[\tilde{v}]))= \frac{n-2}{2}\varphi_0 \left(t^2\tilde{v}^2
+\varepsilon^2\right)^{p/2} \le C_0 \quad \mbox{on }
{\mathcal M},\\
\\
-C_1\le \frac{\p \tilde{v}}{\p \nu}-\frac{n-2}{2}h_g
\tilde{v}=-\frac{n-2}{2} h_0 \left(t^2 \tilde{v}^2
+\varepsilon^2\right)^{q/2}\le -\frac{n-2}{2}h_0 \varepsilon^q
\quad \mbox{on } \p {\mathcal M},
\end{array}\right.
$$
for some positive constants $C_0$ and $C_1$. By Lemma \ref{L6.1}
we have $1/C\le \tilde{v}\le C$ and hence Lemma \ref{L5.1} gives
$\|\tilde{v}\|_{C^2({\mathcal M})}<R_0$ for some number $R_0$
independent of $t$. Hence $\|v\|_{C^2({\mathcal M})}< R_0$. Thus, in
view of (\ref{5.7}),  Claim 2 holds with this $R_0$.

\vskip 0.2cm

From Claim 2 it follows that the Leray-Schauder degree
$\deg(I-T_t, \Phi\cap B_{R_0}, 0)$ is well-defined for $0\le t\le 1$
and is independent of $t$. Since $T_0=0$ we have
$$
\deg(I-T_1, \Phi\cap B_{R_0}, 0)=\deg(I, \Phi\cap B_{R_0}, 0)=1
$$

On the other hand, from Claim 1, Lemma \ref{L6.1}, and Lemma
\ref{L5.1} it follows that there exists $R'\ge R_0$ such that
$\|v\|_{C^2({\mathcal M})}<R'$ for any solution $v$ of ($P_t$)
with $1\le t\le \bar{t}$. Therefore $\deg(I-T_t, B_{R'},0)$ is
well-defined for $t\in [1, \bar{t}]$ and is independent of $t$.
Since ($P_{\bar{t}}$) has no solution, we therefore have
$$
\deg(I-T_1, B_{R'}, 0)=\deg(I-T_{\bar{t}}, B_{R'}, 0)=0.
$$
Let $K:=\bar{B}_{R'}\backslash \Phi$ which is closed in
$C^2({\mathcal M})$. If $0\not\in (I-T_1)(K)$, then the excision
property of Leray-Schauder degree implies that
$$
0=\deg(I-T_1, B_{R'}, 0)=\deg(I-T_1, \Phi\cap B_{R'},0)=1.
$$
which is absurd. Therefore $T_1$ has a fixed point $v_\varepsilon$
in $K$, which is also a solution of ($P_1$). Moreover, the
definition of $\Phi$ implies that $ \sup_{\mathcal M}
v_\varepsilon\ge \a_0>0$. It then follows from Claim 1, Lemma
\ref{L5.3}, Lemma \ref{L5.1}, and the result in \cite{LT86} that
$$
\frac{1}{C}\le v_\varepsilon\le C \quad \mbox{and} \quad
\|v_\varepsilon\|_{C^{2,\a}({\mathcal M})}\le C
$$
for some positive constant $C$ independent of $\varepsilon$, where
$\a\in (0, 1)$. This implies that there is a sequence
$\varepsilon_j\searrow 0$ such that $v_{\varepsilon_j}$ converges
in $C^2({\mathcal M})$ to a solution $v$ of (\ref{5.5}).
\end{proof}

\section{\bf Appendix}
\setcounter{equation}{0}

We include here a comparison principle which is repeatedly used in
the proof of Theorem \ref{T1.2}. For a positive function $v\in
C^2({\mathcal M})$ we still use the notation $V[v]$ defined by
(\ref{5.10}).

\begin{Lemma}\label{L6.1}
Assume that $(f, \Ga)$ satisfies (\ref{1.3}) and (\ref{1.5}) and
that $({\mathcal M}, g)$ is a smooth compact Riemannian manifold with
smooth boundary $\p {\mathcal M}$. Let $a\in C^0(\p {\mathcal M})$ and
let $v, \xi\in C^2({\mathcal M}) $ be two positive functions
satisfying
\begin{equation}\label{A.1}
\left\{\begin{array}{lll} f(\lambda(V[v]))\ge
f(\lambda(V[\xi])),\quad \la(V[v]), \la(V[\xi]) \in \Gamma  \quad
\mbox{on } {\mathcal M},\\
\\
\left(\frac{\p}{\p \nu} -a\right) v\le \min\left\{
\left(\frac{\p}{\p \nu} -a\right)\xi, 0\right\} \quad \mbox{on }
\partial {\mathcal M}.
\end{array}\right.
\end{equation}
Then either $v\equiv \xi$ on ${\mathcal M}$ or $v>\xi$ on ${\mathcal M}$.

\end{Lemma}
\begin{proof} We first prove that $v\ge \xi$ on ${\mathcal M}$.
Suppose it is not true then by using the positivity of $v$, we find a number $\beta>1$
such that $\beta v\ge \xi$ on ${\mathcal M}$ and $\beta v(\bar
x)=\xi(\bar x)$ for some $\bar{x}\in {\mathcal M}$. If $\bar x\in
{\mathcal M}\backslash \p {\mathcal M}$, then
$$
\beta v=\xi, \quad \nabla (\beta v)= \nabla \xi \quad \mbox{and}
\quad \nabla^2 (\beta v) \ge \nabla^2 \xi \quad \mbox{at } \bar{x}
$$
and therefore
$$
\lambda (V[\beta v])\le \lambda(V[\xi]) \quad \mbox{at } \bar{x}.
$$
It follows, using (\ref{1.5}) and $\beta>1$, that
$$
f(\lambda (V[\xi])\ge f(\lambda(V[\beta v])) = f(\beta
\lambda(V[v]))>f(\lambda(V[v])) \quad \mbox{at } \bar{x}
$$
which is a contradiction.

If $\bar x\in \partial {\mathcal M}$, then
$$
\beta v= \xi \quad \mbox{and} \quad \frac{\p}{\p \nu}(\beta
v-\xi)\ge 0 \quad \mbox{at } \bar {x}.
$$
One the other hand, using hypothesis in the lemma, the fact
$\beta>1$, and the fact $(\frac{\p}{\p \nu} -a)v\le 0$ on
$\p{\mathcal M}$, we have
$$
\left(\frac{\p}{\p \nu}-a\right)(\beta v) \le \left(\frac{\p}{\p
\nu}-a\right)v\le \left(\frac{\p}{\p \nu}-a\right)\xi \quad
\mbox{on } \p {\mathcal M}.
$$
Since $(\beta v-\xi)(\bar x)=0$, we have $\frac{\p}{\p \nu}(\beta
v-\xi)\le 0$ at $\bar{x}$. Therefore
$$
\frac{\p}{\p \nu}(\beta v-\xi)(\bar x)= 0.
$$
Using the Hopf Lemma, and the strong maximum principle, we must
have $\beta v-\xi \equiv 0$. This implies
$f(\la(V[\xi]))>f(\la(V[v]))$ on ${\mathcal M}$ as before. We get
a contradiction again.

Therefore $v\ge \xi$ on ${\mathcal M}$. If $v>\xi$ on ${\mathcal
M}$, we are done. Otherwise, there exists $\bar{x}\in {\mathcal
M}$ such that $v(\bar{x})=\xi(\bar{x})$. It then follows from the
boundary condition in (\ref{A.1}) that $\frac{\p}{\p
\nu}(v-\xi)\le 0$ at $\bar{x}$. The Hopf Lemma implies that this
can not occur unless $v\equiv\xi$ on ${\mathcal M}$. If
$\bar{x}\in {\mathcal M}\backslash \p {\mathcal M}$, the strong
maximum principle implies $v\equiv \xi$ on $\p {\mathcal M}$. The
proof is thus complete.
\end{proof}

\end{document}